\documentclass[11pt]{article}
\usepackage{amssymb, amsmath, latexsym, amsbsy}
\usepackage{amsfonts}
\usepackage{dsfont}
\usepackage[T1]{fontenc}
\usepackage[all]{xy}
\usepackage{amscd}
\usepackage{graphicx}
\usepackage{latexsym}

\usepackage{pxfonts}
\usepackage{enumerate}
\usepackage{color}
\usepackage[colorlinks]{hyperref}
\usepackage{fancybox}
\usepackage{fancyhdr}
\usepackage{lipsum}

\newtheorem{definition}{Definition}
\newtheorem{theorem}{Theorem}
\newtheorem{lemma}{Lemma}
\newtheorem{corollary}{Corollary}
\newtheorem{remark}{Remark}
\newtheorem{proposition}{Proposition}
\newenvironment{prof}[1][proof]{\textbf{#1:} }{\ \rule{0.5em}{0.5em}}
\date{\empty}

\date{}
 \title{The prescribed Ricci curvature problem on $5$-dimensional nilpotent Lie groups.}

\begin{document}

 \maketitle
\begin{center}
\author{ \textbf{M. L. Foka}$^{1}$,\quad \textbf{R. P.  Nimpa}$^{2}$,\quad \textbf{S.J. Mbatakou}$^{3}$, \quad \textbf{M. B. N.  Djiadeu}$^{4}$,\quad \textbf{T.B. Bouetou}$^{5}$,\\
\small{e-mail: $\textbf{1}.$ lanndrymarius@gmail.com,\quad $\textbf{2}.$ romain.nimpa@facsciences-uy1.cm, \\ $\textbf{3}.$ salomon-joseph.mbatakou@facsciences-uy1.cm,  \quad
	$\textbf{4}.$ michel.djiadeu@facsciences-uy1.cm, \\
	$\textbf{5}.$ tbouetou@gmail.com,  \\
 University of Yaounde 1, Faculty of Science, Department of
Mathematics,  P.O. Box 812, Yaounde, Republic of Cameroon.}}
\end{center}


\begin{abstract}
	
In this paper, using the Milnor-type theorem  technique, we provide  on each nilpotent five dimensional Lie group, some global existence result of a pair $(g,c)$ consisting of a left-invariant Riemannian metric $g$ and a positive constant $c$ such that $Ric(g)= c T$, where $Ric(g)$ is Ricci curvature of $g$ and  $T$ a given left-invariant symmetric $(0, 2)$-tensor field.
	
  
\end{abstract}

\textbf{keywords}: Lie algebra, Lie group, Milnor-type theorem, Ricci curvature.


\section{Introduction }
An important problem in Riemannian geometry is finding a Riemannian metric $g$ to satisfied the prescribed Ricci curvature equation \begin{equation} \label{riccia}
		Ric(g) = T,
	\end{equation} for some fixed symetric $(0,2)-$tensor field $T$ on a manifold $M$. Some progress on this question concerning the local and global solvability was made in \cite{pu1}, \cite{pu} and \cite{apu}. R. Hamilton in \cite{R.H} and D. DeTurck in \cite{turk}  showed that instead of trying to prove the existence of a metric $g$ solution of \eqref{riccia} , one should look for a metric $g$ and a positive constant $c$ such that \begin{equation} \label{riccib}
Ric(g) = c T,
\end{equation}

Several mathematicians have investigated and obtained some global existence results of \eqref{riccib} under the assumption that the metric and the tensor $T$ are invariant under the transitive action of a Lie group G on M.
	
For non-compact homogeneous spaces, in particular for the connected unimodular Lie groups, the authors of \cite{mil} and \cite{djia}	 observed that
	 the solvabilty of \eqref{riccib} is related to finding the signatures of Ricci curvature of left invariant Riemannian metrics or to study the space of metrics with prescribed Ricci operator.
	 
	  In \cite{{butt}}, Buttsworth settled the question of the solvability of equation \eqref{riccib} for left invariant Riemannian metrics in the case of $3$-dimensional unimodular Lie groups. His results suggested the study of the higher  dimension. Then on non compact homogenous spaces, precisely on five dimensional nilpotent Lie groups, we give in this work the global existence theorems  of  the equation \eqref{riccib}. In order to obtain a suitable  system of polynomials equations, we consider in the sequel the equivalent equation: 
  \begin{equation} \label{riccic}
  	Ric(g) = t^2 T \quad \quad with \: t\in \mathbb{R}^*.
\end{equation}

In section $2$ we gives some preliminaries on 
  on Milnor-type theorem.
In section $3$, we determine explicitly Set of representatives of the moduli space for $5$-dimensional nilpotent Lie algebras. In section $4$, we apply Milnor type theorem to study  all inner products on each $5$-dimensional nilpotent Lie algebra. In section $5$, we compute the Ricci curvature  of left invariant Riemannian metrics for each $5$-dimensional  simply connected nilpotent Lie group. 
  In section $6$, we give a global existence theorem for left
 invariant Riemannian metrics with precribed Ricci curvature on each $5$-dimensional connected nilpotent Lie
 group.

\section{Preliminary on Milnor-type theorems}

	The purpose of this section is to describe the  Milnor-type theorem which is a tool based on the Moduli space of left-invariant Riemannian metrics to provide  the set of all left-invariant Riemannian metrics on a Lie group.  See \cite{{hash},{taka},{koda}} for more details.\\
	
	For a Lie group $G$ of dimension $n$ and  it's Lie algebra $\mathfrak{g}$, the set of left-invariant Riemannnian metrics on $G$ is naturally identified with	
$$\widetilde{\mathfrak{M}}:=\{\langle,\rangle; \: \text{ an inner product on}\;\mathfrak{g}\}.$$
Let's identify $\mathfrak{g}\cong\mathbb{R}^n$ as vector space. The application $\psi:GL_n(\mathbb{R})\times\widetilde{\mathfrak{M}}\longrightarrow \widetilde{\mathfrak{M}}$ defined by
$$\forall g\in GL_n(\mathbb{R}), \:\; \langle,\rangle\ \in \widetilde{\mathfrak{M}},\: \quad\psi(g,\langle,\rangle)= \langle g^{-1}(\cdot),g^{-1}(\cdot)\rangle $$
 induce a transitive action of $GL_n(\mathbb{R})$ on $\widetilde{\mathfrak{M}}$.
 Therefore we have the following identification
 $$\widetilde{\mathfrak{M}}\cong GL_n(\mathbb{R})/O(n).$$
 Let us denote respectively by 
 $$Aut(\mathfrak{g}):= \{\varphi \in GL_n(\mathbb{R})\mid\varphi([\cdot,\cdot])=[\varphi(\cdot),\varphi(\cdot)] \}$$and
 $$\mathbb{R}^\times :=\{c\cdot id:\mathfrak{g}\longrightarrow\mathfrak{g}\mid c\in \mathbb{R}^*\}.$$ 
 The automorphims group of $\mathfrak{g}$ and the scalar group. Note that the group $\mathbb{R}^\times$$Aut(\mathfrak{g})$ naturally acts on  $\widetilde{\mathfrak{M}}$ as subgroup of $GL_n(\mathbb{R})$. 
 
 \begin{definition}
 	The moduli space of left-invariant Riemannian metrics on $G$ is the orbit space of the action of $\mathbb{R}^\times$$Aut(\mathfrak{g})$  on $\widetilde{\mathfrak{M}}$. It is denoted by 
$$\mathfrak{P}\mathfrak{M}\cong \mathbb{R}^\times Aut(\mathfrak{g})/\widetilde{\mathfrak{M}}.$$ \end{definition}
  In what follows, the orbit space of $\mathbb{R}^\times$$Aut(\mathfrak{g})$ through any inner product $\langle,\rangle\in\widetilde{\mathfrak{M}}$ and the double coset of $g\in GL_n(\mathbb{R})$ are respectively denoted by
  $$[\langle,\rangle]:=\mathbb{R}^\times Aut(\mathfrak{g})\cdot\langle,\rangle:=\{\varphi\cdot\langle,\rangle\mid \varphi\in \mathbb{R}^\times Aut(\mathfrak{g})\}$$ and
  $$[[g]]:=\mathbb{R}^\times Aut(\mathfrak{g})gO(n):=\{\varphi gh\mid \varphi\in \mathbb{R}^\times Aut(\mathfrak{g}),h\in O(n)\}.$$
  \begin{definition}
  	A set of representatives of $\mathfrak{P}\mathfrak{M}$ is a subset $U\subset GL_n(\mathbb{R})$ such that $$\mathfrak{P}\mathfrak{M}:=\{[h\cdot\langle,\rangle_0]\mid h\in U\}.$$
  	Where $\langle,\rangle_0$ denote the canonical inner product on $\mathfrak{g}\cong\mathbb{R}^n.$
  \end{definition}  
The  lemma below  gives a criteria for a subset $U\subset GL_n(\mathbb{R})$ to be a set of representatives of $\mathfrak{P}\mathfrak{M}$.
 \begin{lemma}\cite{taka}
 	A subset $U\subset GL_n(\mathbb{R})$ is a set of representatives of $\mathfrak{P}\mathfrak{M}$ if and only if for every $g\in GL_n(\mathbb{R})$ there exist $h\in U$ such that $h\in[[g]]$.
 \end{lemma}
\begin{theorem}\cite{taka}\label{theo01} 
	Let U be a set of representatives of $\mathfrak{PM}$. Then, for every inner
	product $\langle,\rangle$ on $\mathfrak{g}$, There exist $h\in U$, $\varphi \in Aut(\mathfrak{g})$ and.
	$k>0$ such that $\{\varphi he_1,\cdots,\varphi he_n\}$ is an orthonormal basis of $\mathfrak{g}$ whith respect to $k\langle,\rangle$.\end{theorem}
\section{Set of representatives of $\mathfrak{PM}$ for $5$-dimensional nilpotent Lie algebras}
In this section, we give a set of representatives of $\mathfrak{PM}$ for $\mathfrak{g}$, when $\mathfrak{g}$ represents any of the Lie algebras contained in the Table below.\\
\begin{center}\label{tab1}
	Table 
\end{center}
\begin{center}
	\begin{tabular}{lll}
		\hline
		Lie algebra &  &Nonzero commutation relation   \\ 
		\hline$5A_1$ & & none  \\	
		$A_{5,4}$ & &$[e_1,e_4]=e_5$,\; $[e_2,e_3]=e_5$   \\
		$ A_{3,1}\oplus 2A_{1}$& & $[e_1,e_2]=e_5$\\
		$A_{4,1}\oplus A_{1}$& &$[e_1,e_2]=e_3$,\; $[e_1,e_3]=e_5$\\
		$ A_{5,6}$& &$[e_1,e_2]=-e_3$,\; $[e_1,e_3]=e_4$,\; $[e_1,e_4]=e_5$,\; $[e_2,e_3]=e_5$\\
		$A_{5,5}$ & &$[e_1,e_2]=e_4$,\; $[e_1,e_3]=e_5$\; $[e_2,e_4]=e_5$ \\
		$A_{5,3}$ & &$[e_1,e_2]=e_3$,\; $[e_1,e_3]=e_4$\; $[e_2,e_3]=e_5$ \\
		$ A_{5,1}$&  &$[e_1,e_2]=e_4$,\; $[e_1,e_3]=e_5$\\
		$ A_{5,2}$&  &$[e_1,e_2]=e_3$,\; $[e_1,e_3]=e_4$\; $[e_1,e_4]=e_5$\\
		\hline
	\end{tabular}
\end{center}
\begin{remark}
	The above table was reorganised from the one of \cite{nikitenko} in order to obtain a nice matrix expression of a derivation.
\end{remark}
In order to find $\mathfrak{PM}$, one needs a matrix expression of the elements of $Der(\mathfrak{g})$ for each of these Lie algebras.
\begin{equation}\label{eq01}
Der(\mathfrak{g}) := \{D\in \mathfrak{g}l(\mathfrak{g}) \,|\, D[\cdotp, \cdotp] = [D(\cdotp), \cdotp] + [\cdotp, D(\cdotp)]\}
\end{equation} 
The following results are obtained using \eqref{eq01}.
\begin{lemma}\label{lem 001}
	The matrix expression of elements of $Der(\mathfrak{g})$ with respect to the basis
	$\{e_1, e_2, e_3, e_4, e_5\}$ in the previous Table  is given as follows: 
	\begin{itemize}
		
		\item [(1)] For $\mathfrak{g}=5A_{1}$
		\begin{equation}\label{equation2}
		Der(\mathfrak{g})=\left\{{\renewcommand{\arraystretch}{0.4}\begin{pmatrix}
			a_{11}&a_{12}&a_{13}&a_{14}&a_{15}\\
			a_{21}&a_{22}&a_{23}&a_{24}&a_{25}\\
			a_{31}&a_{32}&a_{33}&a_{34}&a_{35}\\
			a_{41}&a_{42}&a_{43}&a_{44}&a_{45}\\
			a_{51}&a_{52}&a_{53}&a_{54}&a_{55}
			\end{pmatrix}}\; ; a_{ij}\in \mathbb{R}  \right\}. \\\end{equation}

		\item [(2)] For $\mathfrak{g}=A_{5,4}$
		\begin{equation}\label{eq02}
		Der(\mathfrak{g})=\left\{{\renewcommand{\arraystretch}{0.4}\begin{pmatrix}
			a_{11}&a_{12}&a_{13}&a_{14}&0\\
			-a_{43}&a_{22}&a_{23}&a_{13}&0\\
			a_{31}&a_{32}&a_{33}&-a_{12}&0\\
			a_{41}&a_{31}&a_{43}&a_{44}&0\\
			a_{51}&a_{52}&a_{53}&a_{54}&a_{55}
			\end{pmatrix}}\; ; a_{ij}\in \mathbb{R}, a_{55}=a_{11}+a_{44}=a_{22}+a_{33}  \right\}. \\\end{equation}
		\item [(3)] For $\mathfrak{g}= A_{3,1}\oplus 2A_{1}$
		\begin{equation}\label{eq03}
		Der(\mathfrak{g})=\left\{{\renewcommand{\arraystretch}{0.4}\begin{pmatrix} a_{11}&a_{12}&0&0&0\\
			a_{21}&a_{22}&0&0&0\\
			a_{31}&a_{32}&a_{33}&a_{34}&0\\
			a_{41}&a_{42}&a_{43}&a_{44}&0\\
			a_{51}&a_{52}&a_{53}&a_{54}&a_{55}
			\end{pmatrix}}\; ; a_{ij}\in \mathbb{R}, a_{55}=a_{11}+a_{22} \right\}. \\\end{equation}
		\item [(4)] For $\mathfrak{g}= A_{4,1}\oplus A_{1}$
		\begin{equation}\label{eq04}
		Der(\mathfrak{g})=\left\{{\renewcommand{\arraystretch}{0.4}\begin{pmatrix} a_{11}&0&0&0&0\\
			a_{21}&a_{22}&0&0&0\\
			a_{31}&a_{32}&a_{33}&0&0\\
			a_{41}&a_{42}&0&a_{44}&0\\
			a_{51}&a_{52}&a_{32}&a_{54}&a_{55}
			\end{pmatrix}}\; ; a_{ij}\in \mathbb{R}, a_{33}=a_{11}+a_{22}, a_{55}=2a_{11}+a_{22} \right\}. \\\end{equation}
		\item [(5)] For $\mathfrak{g}=A_{5,6}$
		\begin{equation}\label{eq05}
		Der(\mathfrak{g})=\left\{{\renewcommand{\arraystretch}{0.4}\begin{pmatrix} a_{11}&0&0&0&0\\
			a_{21}&2a_{11}&0&0&0\\
			a_{31}&a_{32}&3a_{11}&0&0\\
			a_{41}&a_{42}&-a_{32}&4a_{11}&0\\
			a_{51}&a_{52}&a_{53}&a_{54}&5a_{11}
			\end{pmatrix}}\; ; a_{ij}\in \mathbb{R}, a_{21}=a_{54}+a_{32}, a_{53}=a_{31}-a_{42} \right\}. \\\end{equation}
		\item [(6)] For $\mathfrak{g}=A_{5,5}$
		\begin{equation}\label{eq06}
		Der(\mathfrak{g})=\left\{{\renewcommand{\arraystretch}{0.4}\begin{pmatrix}
			a_{11}&a_{12}&0&0&0\\
			a_{21}&a_{22}&0&0&0\\
			a_{31}&a_{32}&2a_{22}&0&0\\
			a_{41}&a_{42}&-a_{12}&a_{44}&0\\
			a_{51}&a_{52}&a_{53}&a_{54}&a_{55}
			\end{pmatrix}}\;a_{ij}\in \mathbb{R} ;a_{44}=a_{11}+a_{22},a_{55}=a_{11}+2a_{22},a_{54}=a_{32}-a_{41} \right\}. \\\end{equation}
		\item [(7)] For $\mathfrak{g}=A_{5,3}$
		\begin{equation}\label{eq07}
		Der(\mathfrak{g})=\left\{{\renewcommand{\arraystretch}{0.4}\begin{pmatrix}
			a_{11}&a_{12}&0&0&0\\
			a_{21}&a_{22}&0&0&0\\
			a_{31}&a_{32}&a_{33}&0&0\\
			a_{41}&a_{42}&a_{32}&a_{44}&a_{12}\\
			a_{51}&a_{52}&-a_{31}&a_{21}&a_{55}
			\end{pmatrix}}\; ;a_{33}=a_{11}+a_{22}, a_{44}=2a_{11}+a_{22},a_{55}=a_{11}+2a_{22} \right\}. \\\end{equation}
		\item [(8)] For $\mathfrak{g}=A_{5,1}$
		\begin{equation}\label{eq08}
		Der(\mathfrak{g})=\left\{{\renewcommand{\arraystretch}{0.4}\begin{pmatrix} a_{11}&0&0&0&0\\
			a_{21}&a_{22}&a_{23}&0&0\\
			a_{31}&a_{32}&a_{33}&0&0\\
			a_{41}&a_{42}&a_{43}&a_{44}&a_{23}\\
			a_{51}&a_{52}&a_{53}&a_{32}&a_{55}
			\end{pmatrix}}\; ; a_{ij}\in \mathbb{R}, a_{55}=a_{11}+a_{33}, a_{44}=a_{11}+a_{22} \right\}. \\\end{equation}
		\item [(9)] For $\mathfrak{g}=A_{5,2}$
		\begin{equation}\label{eq09}
		Der(\mathfrak{g})=\left\{{\renewcommand{\arraystretch}{0.4}\begin{pmatrix}
			a_{11}&0&0&0&0\\
			a_{21}&a_{22}&0&0&0\\
			a_{31}&a_{43}&a_{33}&0&0\\
			a_{41}&a_{42}&a_{43}&a_{44}&0\\
			a_{51}&a_{52}&a_{42}&a_{43}&a_{55}
			\end{pmatrix}}\; ;a_{33}=a_{11}+a_{22}, a_{44}=2a_{11}+a_{22},a_{55}=a_{11}+2a_{22} \right\}. \\\end{equation}
		
	\end{itemize}
\end{lemma}
In the following subsections, we give an explicit expression of $\mathfrak{PM}$ for every lie algebra $\mathfrak{g}$ contains in Table. We fix a basis
$\{e_1, e_2, e_3, e_4, e_5\}$ of $\mathfrak{g}$ whose brackets relations are given in the Table .
\subsection{Case of $\mathfrak{g}=5A_{1}$}

Let's consider $\mathbb{R}^{\times}Aut(\mathfrak{g})$ for this Lie algebra. From the lemma \ref{lem 001}  one can see that $\mathbb{R}\oplus Der(\mathfrak{g})$ which is the Lie algebra of $\mathbb{R}^{\times}Aut(\mathfrak{g})$ contains 

\begin{equation*}G_1=\left\{\renewcommand{\arraystretch}{0.4}\begin{pmatrix}
0&0&0&0&0\\
a_{21}&0&0&0&0\\
a_{31}&a_{32}&0&0&0\\
a_{41}&a_{42}&a_{43}&0&0\\
a_{51}&a_{52}&a_{53}&a_{54}&0
\end{pmatrix};\;a_{ij}\in\mathbb{R}\right\}.
\end{equation*}
Hence, by exponentiation it, we obtain
\begin{equation}\label{equation 4}
F_1=\left\{{\renewcommand{\arraystretch}{0.4}\begin{pmatrix}
	1&0&0&0&0\\
	a_{21}&1&0&0&0\\
	a_{31}&a_{32}&1&0&0\\
	a_{41}&a_{42}&a_{43}&1&0\\
	a_{51}&a_{52}&a_{53}&a_{54}&1
	\end{pmatrix}}\; ;a_{ij}\in \mathbb{R} \right\}\subseteq \mathbb{R}^{\times}Aut(\mathfrak{g}).
\end{equation}
\begin{lemma}\label{lemma2}
	Let $\mathfrak{g}= 5A_{1}$. Then the following $U_1'$ is a set of representatives of $\mathfrak{PM}$	
	\begin{equation}\label{equation3}
	U_1'=\left\{diag(a,b,c,d,e) ; a,b,c,d,e>0\right\}. \\\end{equation}
\end{lemma}

\begin{prof}
	Take any $g \in GL_5(\mathbb{R})$. We need to obtain $a,b,c,d,e>0$ such that: 
	$diag(a,b,c,d,e) \in [[g]]$.\\
	There exists $k\in O(5)$ such that $gk={\renewcommand{\arraystretch}{0.4}\begin{pmatrix}
		g_{11}&0&0&0&0\\
		g_{21}&g_{22}&0&0&0\\
		g_{31}&g_{32}&g_{33}&0&0\\
		g_{41}&g_{42}&g_{43}&g_{44}&0\\
		g_{51}&g_{52}&g_{53}&g_{54}&g_{55}
		\end{pmatrix}}$	 with $g_{11},g_{22},g_{33},g_{44},g_{55}>0$.\\ Let consider the matrix $\varphi={\renewcommand{\arraystretch}{0.4}\begin{pmatrix}
		1&0&0&0&0\\
		a_{21}&1&0&0&0\\
		a_{31}&a_{32}&1&0&0\\
		a_{41}&a_{42}&a_{43}&1&0\\
		a_{51}&a_{52}&a_{53}&a_{54}&1
		\end{pmatrix}} \in F_1\subseteq \mathbb{R}^\times Aut(\mathfrak{g})$. With
	$$a_{21}=-\frac{g_{21}}{g_{11}},\quad
	a_{32}=-\frac{g_{32}}{g_{22}},\quad
	a_{31}=-\frac{g_{31}+a_{32}g_{21}}{g_{11}},\quad
	a_{41}=-\frac{g_{41}+a_{43}g_{31}+a_{42}g_{21}}{g_{11}},$$
	$$a_{42}=-\frac{g_{42}+a_{43}g_{32}}{g_{22}},\quad
	a_{43}=-\frac{g_{43}}{g_{33}},\quad
	a_{51}=-\frac{g_{51}+a_{54}g_{41}+a_{53}g_{31}+a_{52}g_{21}}{g_{11}},$$
	$$a_{52}=-\frac{g_{52}+a_{54}g_{42}+a_{53}g_{32}}{g_{22}},\quad
	a_{53}=-\frac{g_{53}+a_{54}g_{43}}{g_{33}},\quad
	a_{54}=-\frac{g_{54}}{g_{44}}.$$
	
	 Hence $\varphi gk=dig(g_{11},g_{22},g_{33},g_{44},g_{55})\in[[g]]$. \\	
One can take $a=g_{11},b=g_{22},c=g_{33}, d=g_{44}$ and $e=g_{55}$ to conclude the proof. 
\end{prof}

\begin{proposition}\label{proposition1}Let $\mathfrak{g}=5A_{1}$. Then $U_1=\left\{I_5\right\} $ is a set of representatives of $\mathfrak{PM}$.
\end{proposition}
\begin{prof}
	By Lemma \ref{lemma2}, there exists $a,b,c,d,e>0$ such that $h=diag(a,b,c,d,e)$. \\It follows from  \eqref{equation2} that, $A=diag\Big(\frac{1}{a},\frac{1}{b},\frac{1}{c},\frac{1}{d},\frac{1}{e}\Big)
	\in \mathbb{R}^{\times} Aut(\mathfrak{g})$.\\
	Hence, $Ah=diag(1,1,1,1,1)\in[[g]]$. 	
\end{prof}

\subsection{Case of $\mathfrak{g}=A_{5,4}$}
Let's consider $\mathbb{R}^{\times}Aut(\mathfrak{g})$ for this Lie algebra. From the lemma \ref{lem 001}  one can see that $\mathbb{R}\oplus Der(\mathfrak{g})$ which is the Lie algebra of $\mathbb{R}^{\times}Aut(\mathfrak{g})$ contains 

\begin{equation*}G_2=\left\{\renewcommand{\arraystretch}{0.4}\begin{pmatrix}
0&0&0&0&0\\
-a_{43}&0&0&0&0\\
a_{31}&a_{32}&0&0&0\\
a_{41}&a_{31}&a_{43}&0&0\\
a_{51}&a_{52}&a_{53}&a_{54}&0
\end{pmatrix};\;a_{ij}\in\mathbb{R}\right\}.
\end{equation*}
Hence, by matrix exponentiation, we obtain
\begin{equation}\label{e0004}
F_2=\left\{{\renewcommand{\arraystretch}{0.4}\begin{pmatrix}
	1&0&0&0&0\\
	-a_{43}&1&0&0&0\\
	a_{31}&a_{32}&1&0&0\\
	a_{41}&a_{42}&a_{43}&1&0\\
	a_{51}&a_{52}&a_{53}&a_{54}&1
	\end{pmatrix}}\; ;a_{ij}\in \mathbb{R} \right\}\subseteq \mathbb{R}^{\times}Aut(\mathfrak{g}).
\end{equation}
\begin{lemma}\label{lem00002}
	Let $\mathfrak{g}= A_{5,4}$. Then the following $U'_2$ is a set of representatives of $\mathfrak{PM}$	
	\begin{equation}\label{e0003}
	U'_2=\left\{{\renewcommand{\arraystretch}{0.4}\begin{pmatrix}
		a&0&0&0&0\\
		f&b&0&0&0\\
		0&0&c&0&0\\
		0&0&0&d&0\\
		0&0&0&0&e
		\end{pmatrix}}\; ; a,b,c,d,e>0,\; f\in \mathbb{R}\right\}. \\\end{equation}
\end{lemma}

\begin{prof}
	Take any $g \in GL_5(\mathbb{R})$. Let's find  $a,b,c,d,e>0$, $f\in \mathbb{R}$ such that: 
	\begin{equation*}
	{\renewcommand{\arraystretch}{0.4}\begin{pmatrix}
		a&0&0&0&0\\
		f&b&0&0&0\\
		0&0&c&0&0\\
		0&0&0&d&0\\
		0&0&0&0&e
		\end{pmatrix}} \in [[g]].
	\end{equation*}
	There exists $k\in O(5)$ such that $gk={\renewcommand{\arraystretch}{0.2}\begin{pmatrix}
		g_{11}&0&0&0&0\\
		g_{21}&g_{22}&0&0&0\\
		g_{31}&g_{32}&g_{33}&0&0\\
		g_{41}&g_{42}&g_{43}&g_{44}&0\\
		g_{51}&g_{52}&g_{53}&g_{54}&g_{55}
		\end{pmatrix}}$	 with $g_{11},g_{22},g_{33},g_{44},g_{55}>0$.\\ Let consider the matrix $\varphi={\renewcommand{\arraystretch}{0.4}\begin{pmatrix}
		1&0&0&0&0\\
		-a_{43}&1&0&0&0\\
		a_{31}&a_{32}&1&0&0\\
		a_{41}&a_{42}&a_{43}&1&0\\
		a_{51}&a_{52}&a_{53}&a_{54}&1
		\end{pmatrix}} \in F_2\subseteq \mathbb{R}^\times Aut(\mathfrak{g})$ with
	$$a_{32}=-\frac{g_{32}}{g_{22}}, \quad
	a_{31}=-\frac{g_{31}+a_{32}g_{21}}{g_{11}},\quad
	a_{41}=-\frac{g_{41}+a_{43}g_{31}+a_{42}g_{21}}{g_{11}},\quad
	a_{42}=-\frac{g_{42}+a_{43}g_{32}}{g_{22}},\quad
	a_{43}=-\frac{g_{43}}{g_{33}},$$
	$$a_{51}=-\frac{g_{51}+a_{54}g_{41}+a_{53}g_{31}+a_{52}g_{21}}{g_{11}},\quad
	a_{52}=-\frac{g_{52}+a_{54}g_{42}+a_{53}g_{32}}{g_{22}},\quad
	a_{53}=-\frac{g_{53}+a_{54}g_{43}}{g_{33}},\quad
	a_{54}=-\frac{g_{54}}{g_{44}}.$$

Hense $\varphi gk={\renewcommand{\arraystretch}{0.3}\begin{pmatrix}
		g_{11}&0&0&0&0\\
		-a_{43}g_{11}+g_{21}&g_{22}&0&0&0\\
		0&0&g_{33}&0&0\\
		0&0&0&g_{44}&0\\
		0&0&0&0&g_{55}
		\end{pmatrix}}\in [[g]]$. \\
	
	One can Take: $a=g_{11},b=g_{22},c=g_{33},d=g_{44}$, $e=g_{55}$, $f=-a_{43}g_{11}+g_{21}$ to conclude the proof. 
\end{prof}

\begin{proposition}\label{propo00001}Let $\mathfrak{g}=A_{5,4}$. Then the following $U_2$ is a set of representatives of $\mathfrak{PM}$
	\begin{equation}\label{eq00004}	
	U_2=\left\{{\renewcommand{\arraystretch}{0.4}\begin{pmatrix}
		1&0&0&0&0\\
		a_{21}&1&0&0&0\\
		0&0&1&0&0\\
		0&0&0&a_{44}&0\\
		0&0&0&0&a_{55}
		\end{pmatrix}}\; ; a_{44},a_{55}>0,\; a_{21}\in\mathbb{R}\right\} \\\end{equation}.
\end{proposition}
\begin{prof}
	By Lemma \ref{lem00002}, there exists $a,b,c,d,e>0$, $f\in \mathbb{R}$ such that $h={\renewcommand{\arraystretch}{0.4}\begin{pmatrix}
		a&0&0&0&0\\
		f&b&0&0&0\\
		0&0&c&0&0\\
		0&0&0&d&0\\
		0&0&0&0&e
		\end{pmatrix}} \in [[g]]$. \\It follows from  \eqref{eq02} that, $A=diag\Big(\frac{1}{a},\frac{1}{b},\frac{1}{c},\frac{a}{bc},\frac{1}{bc}\Big)
	\in \mathbb{R}^{\times} Aut(\mathfrak{g})$.\\
	Hence, $Ah={\renewcommand{\arraystretch}{0.4}\begin{pmatrix}
		1&0&0&0&0\\
		\frac{f}{b}&1&0&0&0\\
		0&0&1&0&0\\
		0&0&0&\frac{ad}{bc}&0\\
		0&0&0&0&\frac{e}{bc}\end{pmatrix}}\in[[g]]$. Therefore taking  $a_{21}=\frac{f}{b}$;  $a_{44}=\frac{ad}{bc}$ ; $a_{55}=\frac{e}{bc}$ we complete the proof.

\end{prof}

\subsection{Case of $\mathfrak{g}=A_{3,1}\oplus 2A_1$}
Let consider $\mathbb{R}^{\times}Aut(\mathfrak{g})$ for this Lie algebra. From the lemma \ref{lem 001}  one can see that $\mathbb{R}\oplus Der(\mathfrak{g})$ which is the Lie algebra of $\mathbb{R}^{\times}Aut(\mathfrak{g})$ contains 
\begin{equation*}
G_3=\left\{{\renewcommand{\arraystretch}{0.4}\begin{pmatrix}
	0&0&0&0&0\\
	a_{21}&0&0&0&0\\
	a_{31}&a_{32}&0&0&0\\
	a_{41}&a_{42}&a_{43}&0&0\\
	a_{51}&a_{52}&a_{53}&a_{54}&0
	\end{pmatrix}}\; ;a_{ij}\in \mathbb{R}\right\}. \end{equation*}
Hence, by matrix exponentiation , we obtain

\begin{equation}\label{e002}
F_3=\left\{{\renewcommand{\arraystretch}{0.4}\begin{pmatrix}
	1&0&0&0&0\\
	a_{21}&1&0&0&0\\
	a_{31}&a_{32}&1&0&0\\
	a_{41}&a_{42}&a_{43}&1&0\\
	a_{51}&a_{52}&a_{53}&a_{54}&1
	\end{pmatrix}}\; ;a_{ij}\in \mathbb{R} \right\}\subset \mathbb{R}^{\times}Aut(\mathfrak{g}).\\\end{equation}
\begin{lemma}\label{lem003}
	Let $\mathfrak{g}=A_{3,1}\oplus 2A_1$. Then the following $U'_3$ is a set of representatives of $\mathfrak{PM}$	
	\begin{equation}\label{e003}
	U'_3=\left\{diag(a,b,c,d,e)
	\; ; a,b,c,d,e>0 \right\}. \\\end{equation}	
\end{lemma}
\begin{prof}
	Take any $g \in GL_5(\mathbb{R})$. Let's find $a,b,c,d,e>0$ such that $diag(a,b,c,d,e)\in [[g]]$.\\
	
	There exists $k\in O(5)$ such that $gk={\renewcommand{\arraystretch}{0.4}\begin{pmatrix}
		g_{11}&0&0&0&0\\
		g_{21}&g_{22}&0&0&0\\
		g_{31}&g_{32}&g_{33}&0&0\\
		g_{41}&g_{42}&g_{43}&g_{44}&0\\
		g_{51}&g_{52}&g_{53}&g_{54}&g_{55}
		\end{pmatrix}}$	 with $g_{11},g_{22},g_{33},g_{44},g_{55}>0$.\\ Let consider the matrix $\varphi={\renewcommand{\arraystretch}{0.4}\begin{pmatrix}
		1&0&0&0&0\\
		a_{21}&1&0&0&0\\
		a_{31}&a_{32}&1&0&0\\
		a_{41}&a_{42}&a_{43}&1&0\\
		a_{51}&a_{52}&a_{53}&a_{54}&1
		\end{pmatrix}}\in F_3\subseteq \mathbb{R}^\times Aut(\mathfrak{g})$ with
	$$a_{21}=-\frac{g_{21}}{g_{11}}, \quad
	a_{31}=-\frac{g_{31}+a_{32}g_{21}}{g_{11}}, \quad
	a_{32}=-\frac{g_{32}}{g_{22}}, \quad
	a_{41}=-\frac{g_{41}+a_{43}g_{31}+a_{42}g_{21}}{g_{11}},$$
	$$a_{42}=-\frac{g_{42}+a_{43}g_{32}}{g_{22}},\quad
	a_{43}=-\frac{g_{43}}{g_{33}},\quad
	a_{51}=-\frac{g_{51}+a_{54}g_{41}+a_{53}g_{31}+a_{52}g_{21}}{g_{11}},\quad
	a_{52}=-\frac{g_{52}+a_{54}g_{42}+a_{53}g_{32}}{g_{22}},$$
	$$a_{53}=-\frac{g_{53}+a_{54}g_{43}}{g_{33}},\quad
	a_{54}=-\frac{g_{54}}{g_{44}}.$$
	
	Hence $\varphi gk=diag(g_{11},g_{22},g_{33},g_{44},g_{55})\in[[g]]$.\\ 
	Take $a=g_{11},b=g_{22},c=g_{33},d=g_{44}$, $e=g_{55}$ to conclude the proof.
\end{prof}
\begin{proposition}\label{propo002}Let $\mathfrak{g}=A_{3,1}\oplus 2A_1$. Then the following $U_1$ is a set of representatives of $\mathfrak{PM}$
	\begin{equation}\label{eq003}	
	U_2=\left\{
	diag(1,1,1,1,a_{55})\; ; a_{55}>0\right\}. \\\end{equation}
\end{proposition}
\begin{prof}
	By Lemma \ref{lem003}, there exists $a,b,c,d,e>0$ such that $h=diag(a,b,c,d,e)
	\in [[g]]$. \\It follows from  \eqref{eq03} that, $A=diag\Big(\frac{1}{a},\frac{1}{b},\frac{1}{c},\frac{1}{d},\frac{1}{ab}\Big)
	\in \mathbb{R}^{\times} Aut(\mathfrak{g})$.\\
	Hence, $Ah=diag\Big(1,1,1,1,\frac{e}{ab}\Big)\in[[g]]$
	. Therefore, taking , $a_{55}=\frac{e}{ab}$ we complete the proof.
	
\end{prof}

\subsection{Case of $\mathfrak{g}=A_{4,1}\oplus A_{1}$}
Let's consider $\mathbb{R}^{\times}Aut(\mathfrak{g})$ for this Lie algebra. From the lemma \ref{lem 001}  one can see that $\mathbb{R}\oplus Der(\mathfrak{g})$ which is the Lie algebra of $\mathbb{R}^{\times}Aut(\mathfrak{g})$ contains 
\begin{equation*}
G_4=\left\{{\renewcommand{\arraystretch}{0.04}\begin{pmatrix}
	0&0&0&0&0\\
	a_{21}&0&0&0&0\\
	a_{31}&a_{32}&0&0&0\\
	a_{41}&a_{42}&0&0&0\\
	a_{51}&a_{52}&a_{32}&a_{54}&0
	\end{pmatrix}}\; ;a_{ij}\in \mathbb{R}\right\}. \end{equation*}
Hence, by matrix exponentiation, we obtain

\begin{equation}\label{e003}
F_4=\left\{{\renewcommand{\arraystretch}{0.4}\begin{pmatrix}
	1&0&0&0&0\\
	a_{21}&1&0&0&0\\
	a_{31}&a_{32}&1&0&0\\
	a_{41}&a_{42}&0&1&0\\
	a_{51}&a_{52}&a_{32}&a_{54}&1
	\end{pmatrix}}\; ;a_{ij}\in \mathbb{R} \right\}\subset \mathbb{R}^{\times}Aut(\mathfrak{g}). \\\end{equation}
\begin{lemma}\label{lem004}
	Let $\mathfrak{g}= A_{4,1}\oplus A_{1}$. Then the following $U'_4$ is a set of representatives of $\mathfrak{PM}$	
	\begin{equation}\label{e003}
	U'_4=\left\{{\renewcommand{\arraystretch}{0.4}\begin{pmatrix}
		a&0&0&0&0\\
		0&b&0&0&0\\
		0&0&c&0&0\\
		0&0&f&d&0\\
		0&0&l&0&e
		\end{pmatrix}}\; ; a,b,c,d,e>0,\; l,f\in \mathbb{R}\right\}. \\\end{equation}
\end{lemma}
\begin{prof}
	Take any $g \in GL_5(\mathbb{R})$. Let's find  $a,b,c,d,e>0$, $f, l\in \mathbb{R}$ such that: 
	\begin{equation*}
	{\renewcommand{\arraystretch}{0.4}\begin{pmatrix}
		a&0&0&0&0\\
		0&b&0&0&0\\
		0&0&c&0&0\\
		0&0&f&d&0\\
		0&0&l&0&e
		\end{pmatrix}} \in [[g]].
	\end{equation*}
	There exists $k\in O(5)$ such that $gk={\renewcommand{\arraystretch}{0.4}\begin{pmatrix}
		g_{11}&0&0&0&0\\
		g_{21}&g_{22}&0&0&0\\
		g_{31}&g_{32}&g_{33}&0&0\\
		g_{41}&g_{42}&g_{43}&g_{44}&0\\
		g_{51}&g_{52}&g_{53}&g_{54}&g_{55}
		\end{pmatrix}}$	 with $g_{11},g_{22},g_{33},g_{44},g_{55}>0$.\\ Let consider the matrix $\varphi={\renewcommand{\arraystretch}{0.4}\begin{pmatrix}
		1&0&0&0&0\\
		a_{21}&1&0&0&0\\
		a_{31}&a_{32}&1&0&0\\
		a_{41}&a_{42}&0&1&0\\
		a_{51}&a_{52}&a_{32}&a_{54}&1
		\end{pmatrix}} \in F_4\subseteq \mathbb{R}^\times Aut(\mathfrak{g})$ with
	$$a_{21}=-\frac{g_{21}}{g_{11}}, \quad
	a_{31}=-\frac{g_{31}+a_{32}g_{21}}{g_{11}}, \quad
	a_{32}=-\frac{g_{32}}{g_{22}}, \quad
	a_{41}=-\frac{g_{41}+a_{42}g_{21}}{g_{11}},\quad
	a_{42}=-\frac{g_{42}}{g_{22}},$$
	$$a_{51}=-\frac{g_{51}+a_{54}g_{41}+a_{32}g_{31}+a_{52}g_{21}}{g_{11}},\quad
	a_{52}=-\frac{g_{52}+a_{32}g_{32}+a_{54}g_{42}}{g_{22}},\quad
	a_{54}=-\frac{g_{54}}{g_{44}}.$$

Hence $\varphi gk={\renewcommand{\arraystretch}{0.4}\begin{pmatrix}
		g_{11}&0&0&0&0\\
		0&g_{22}&0&0&0\\
		0&0&g_{33}&0&0\\
		0&0&g_{43}&g_{44}&0\\
		0&0&-\frac{g_{33}g_{32}}{g_{22}}&0&g_{55}
		\end{pmatrix}}\in[[g]]$. \\ One can Take $a=g_{11},b=g_{22},c=g_{33},d=g_{44}$, $e=g_{55}$, $f=g_{43}$ $l=-\frac{g_{33}g_{32}}{g_{22}}$ to conclude the proof.
\end{prof}
\begin{lemma}\label{propo0003}Let $\mathfrak{g}= A_{4,1}\oplus A_{1}$. Then the following $U"_4$ is a set of representatives of $\mathfrak{PM}$
	\begin{equation}\label{eq004}	
	U"_4=\left\{{\renewcommand{\arraystretch}{0.4}\begin{pmatrix}
		1&0&0&0&0\\
		0&1&0&0&0\\
		0&0&a_{33}&0&0\\
		0&0&a_{43}&1&0\\
		0&0&a_{53}&0&a_{55}
		\end{pmatrix}}\; ; a_{44},a_{55}>0,\; a_{53},a_{43}\in\mathbb{R}\right\}. \\\end{equation}
\end{lemma}
\begin{prof}
	By Lemma \ref{lem004}, there exists $a,b,c,d,e>0$, $f,l\in \mathbb{R}$ such that $h={\renewcommand{\arraystretch}{0.4}\begin{pmatrix}
		a&0&0&0&0\\
		0&b&0&0&0\\
		0&0&c&0&0\\
		0&0&f&d&0\\
		0&0&l&0&e
		\end{pmatrix}} \in [[g]]$. $\forall g\in GL_n(\mathbb{R})$.\\It follows from  \eqref{eq04} that,$A=diag\Big(\frac{1}{a},\frac{1}{b},\frac{1}{ab},\frac{1}{d},\frac{1}{a^2b}\Big)
	\in \mathbb{R}^{\times} Aut(\mathfrak{g})$.\\
	Hence, $Ah={\renewcommand{\arraystretch}{0.4}\begin{pmatrix}
		1&0&0&0&0\\
		0&1&0&0&0\\
		0&0&\frac{c}{ab}&0&0\\
		0&0&\frac{f}{d}&1&0\\
		0&0&\frac{l}{a^2b}&0&\frac{e}{a^2b}
		\end{pmatrix}}\in[[g]]$. Therefore taking ,  $a_{43}=\frac{f}{d}$; $a_{53}=\frac{l}{a^2b}$; $a_{33}=\frac{c}{ab}$ ; $a_{55}=\frac{e}{a^2b}$ we complete the proof.
\end{prof}\\

\begin{remark} When $\mathfrak{g}=A_{4,1}\oplus A_{1}$ we have two cases. Indeed, we have two set of representatives with the same number of parameters. This is due to the fact that every element of the set of representatives $U_3$ given by the previous lemma is entirely defined by the data of $a_{33},a_{55}>0$ and $a_{43},a_{53}\in\mathbb{R}$. So for a given $h\in U_4$, $a_{43}=0$ or $a_{43}\neq0$. The following proposition clearly illustrates the situation.\end{remark}
\begin{proposition}\label{propo003}
	Let $\mathfrak{g}= A_{4,1}\oplus A_{1}$. Then the following $U^1_4$ or $U^2_4$ is a set of representatives of $\mathfrak{PM}$
	\begin{itemize}
		
		\item[1)] 	\begin{equation}\label{eq004}	
		U^1_4=\left\{{\renewcommand{\arraystretch}{0.4}\begin{pmatrix}
			1&0&0&0&0\\
			0&1&0&0&0\\
			0&0&a&0&0\\
			0&0&0&1&0\\
			0&0&c&0&b
			\end{pmatrix}}\; ; a_{44},a_{55}>0,\; a_{53}\in\mathbb{R}\right\} \\\end{equation}
		\item[2)] 	\begin{equation}\label{eq004}	
		U^2_4=\left\{{\renewcommand{\arraystretch}{0.4}\begin{pmatrix}
			1&0&0&0&0\\
			0&1&0&0&0\\
			0&0&a&0&0\\
			0&0&c&1&0\\
			0&0&0&0&b
			\end{pmatrix}}\; ; a_{44},a_{55}>0,\; a_{53}\in\mathbb{R}\right\}. \\\end{equation}
	\end{itemize}
\end{proposition}
\begin{prof}
	By Lemma \ref{propo0003}, there exists $a_{33},a_{55}>0$, $a_{43},a_{53}\in \mathbb{R}$ such that $h={\renewcommand{\arraystretch}{0.4}\begin{pmatrix}
		1&0&0&0&0\\
		0&1&0&0&0\\
		0&0&a_{33}&0&0\\
		0&0&a_{43}&1&0\\
		0&0&a_{53}&0&a_{55}
		\end{pmatrix}} \in [[g]]$. $\forall g\in GL_n(\mathbb{R})$\\
	\begin{enumerate}
		\item [1)] \underline{First case}\\
		If $a_{43}=0$ then $U^1_4$ is a set of representatives of $\mathfrak{PM}$
		\item [2)]\underline{Second case}\\
		If $a_{43}\neq0$ then it follows from  \eqref{eq04} that,$A={\renewcommand{\arraystretch}{0.4}\begin{pmatrix}
			1&0&0&0&0\\
			0&1&0&0&0\\
			0&0&1&0&0\\
			0&0&0&1&0\\
			0&0&0&-\frac{a_{53}}{a_{43}}&1
			\end{pmatrix}} \in \mathbb{R}^{\times} Aut(\mathfrak{g})$.\\
		Hence, $Ah={\renewcommand{\arraystretch}{0.4}\begin{pmatrix}
			1&0&0&0&0\\
			0&1&0&0&0\\
			0&0&a_{33}&0&0\\
			0&0&a_{43}&1&0\\
			0&0&0&0&a_{55}
			\end{pmatrix}}\in[[g]]$. Therefore $U^2_4$ is a set of representatives of $\mathfrak{PM}$.	\end{enumerate}
\end{prof}\\ 

\subsection{Case of $\mathfrak{g}=A_{5,6}$}
Let's consider $\mathbb{R}^{\times}Aut(\mathfrak{g})$ for this Lie algebra. From the lemma \ref{lem 001}  one can see that $\mathbb{R}\oplus Der(\mathfrak{g})$ which is the Lie algebra of $\mathbb{R}^{\times}Aut(\mathfrak{g})$ contains 
\begin{equation*}
G_5=\left\{{\renewcommand{\arraystretch}{0.4}\begin{pmatrix}
	0&0&0&0&0\\
	a_{54}+a_{32}&0&0&0&0\\
	a_{31}&a_{32}&0&0&0\\
	a_{41}&a_{42}&-a_{32}&0&0\\
	a_{51}&a_{52}&a_{31}-a_{42}&a_{54}&0
	\end{pmatrix}}\; ;a_{ij}\in \mathbb{R}\right\}.
\end{equation*}
Hence, by matrix  exponentiation, we obtain

\begin{equation}\label{e004}
F_5=\left\{{\renewcommand{\arraystretch}{0.4}\begin{pmatrix}
	1&0&0&0&0\\
	a_{54}+a_{32}&1&0&0&0\\
	a_{31}&a_{32}&1&0&0\\
	a_{41}&a_{42}&-a_{32}&1&0\\
	a_{51}&a_{52}&a_{53}&a_{21}-a_{32}&1
	\end{pmatrix}}\; ;a_{ij}\in \mathbb{R} \right\}\subseteq \mathbb{R}^{\times}Aut(\mathfrak{g})
\end{equation}
\begin{lemma}\label{lem005}
	Let $\mathfrak{g}= A_{5,6}$. Then the following $U'_5$ is a set of representatives of $\mathfrak{PM}$	
	\begin{equation}\label{e003}
	U'_5=\left\{{\renewcommand{\arraystretch}{0.4}\begin{pmatrix}
		a&0&0&0&0\\
		f&b&0&0&0\\
		0&0&c&0&0\\
		0&0&l&d&0\\
		0&0&0&0&e
		\end{pmatrix}}\; ; a,b,c,d,e>0,\; l,f\in \mathbb{R}\right\}. \\\end{equation}
\end{lemma}

\begin{prof}
	Take any $g \in GL_5(\mathbb{R})$. Let's find  $a,b,c,d,e>0$, $f, l\in \mathbb{R}$ such that: 
	\begin{equation*}
	{\renewcommand{\arraystretch}{0.4}\begin{pmatrix}
		a&0&0&0&0\\
		f&b&0&0&0\\
		0&0&c&0&0\\
		0&0&l&d&0\\
		0&0&0&0&e
		\end{pmatrix}} \in [[g]].
	\end{equation*}
There exists $k\in O(5)$ such that $gk={\renewcommand{\arraystretch}{0.4}\begin{pmatrix}
		g_{11}&0&0&0&0\\
		g_{21}&g_{22}&0&0&0\\
		g_{31}&g_{32}&g_{33}&0&0\\
		g_{41}&g_{42}&g_{43}&g_{44}&0\\
		g_{51}&g_{52}&g_{53}&g_{54}&g_{55}
		\end{pmatrix}}$	 with $g_{11},g_{22},g_{33},g_{44},g_{55}>0$.\\ Let consider the matrix $\varphi={\renewcommand{\arraystretch}{0.4}\begin{pmatrix}
		1&0&0&0&0\\
		a_{21}&1&0&0&0\\
		a_{31}&a_{32}&1&0&0\\
		a_{41}&a_{42}&-a_{32}&1&0\\
		a_{51}&a_{52}&a_{53}&a_{54}&1
		\end{pmatrix}} \in F_5\subseteq \mathbb{R}^\times Aut(\mathfrak{g})$ with
	$$a_{21}=a_{54}+a_{32}, \quad
	a_{31}=-\frac{g_{31}+a_{32}g_{21}}{g_{11}}, \quad
	a_{32}=-\frac{g_{32}}{g_{22}}, \quad
	a_{41}=-\frac{g_{41}+a_{43}g_{31}+a_{42}g_{21}}{g_{11}},\quad
	a_{42}=-\frac{g_{42}+a_{43}g_{32}}{g_{22}},$$
	$$a_{43}=-a_{32},\quad
	a_{51}=-\frac{g_{51}+a_{54}g_{41}+a_{53}g_{31}+a_{52}g_{21}}{g_{11}},\quad
	a_{52}=-\frac{g_{52}+a_{54}g_{42}+a_{53}g_{32}}{g_{22}},$$
	$$a_{53}=-\frac{g_{53}+a_{54}g_{43}}{g_{33}},\quad
	a_{54}=\frac{g_{54}}{g_{44}}.$$	
	
 Hence $\varphi gk={\renewcommand{\arraystretch}{0.3}\begin{pmatrix}
		g_{11}&0&0&0&0\\
		(a_{54}+a_{32})g_{11}+g_{21}&g_{22}&0&0&0\\
		0&0&g_{33}&0&0\\
		0&0&-a_{32}g_{33}+g_{43}&g_{44}&0\\
		0&0&0&0&g_{55}
		\end{pmatrix}}\in[[g]]$. \\
	
	One can Take: $a=g_{11},b=g_{22},c=g_{33},d=g_{44}$, $e=g_{55}$, $l=-a_{32}g_{33}+g_{43}$ $f=(a_{54}+a_{32})g_{11}+g_{21}$ to conclude the proof. 
\end{prof}

\begin{proposition}\label{propo004}Let $\mathfrak{g}=A_{5,6}$. Then the following $U_5$ is a set of representatives of $\mathfrak{PM}$
	\begin{equation}\label{eq004}	
	U_5=\left\{{\renewcommand{\arraystretch}{0.4}\begin{pmatrix}
		1&0&0&0&0\\
		a_{21}&a_{22}&0&0&0\\
		0&0&a_{33}&0&0\\
		0&0&a_{43}&a_{44}&0\\
		0&0&0&0&a_{55}
		\end{pmatrix}}\; ; a_{22},a_{33},a_{44},a_{55}>0,\; a_{21}, a_{43}\in\mathbb{R}\right\}. \\\end{equation}
\end{proposition}
\begin{prof}
	By Lemma \ref{lem005}, there exists $a,b,c,d,e>0$, $f,l\in \mathbb{R}$ such that $h={\renewcommand{\arraystretch}{0.4}\begin{pmatrix}
		a&0&0&0&0\\
		f&b&0&0&0\\
		0&0&c&0&0\\
		0&0&l&d&0\\
		0&0&0&0&e
		\end{pmatrix}} \in [[g]]$. \\It follows from  \eqref{eq05} that,$A=diag\Big(\frac{1}{a},\frac{1}{a^2},\frac{1}{a^3},\frac{1}{a^4},\frac{1}{a^5}\Big)
	\in \mathbb{R}^{\times} Aut(\mathfrak{g})$.\\
	Hence,$Ah={\renewcommand{\arraystretch}{0.4}\begin{pmatrix}
		1&0&0&0&0\\
		\frac{f}{a^2}&\frac{b}{a^2}&0&0&0\\
		0&0&\frac{c}{a^3}&0&0\\
		0&0&\frac{l}{a^4}&\frac{d}{a^4}&0\\
		0&0&0&0&\frac{e}{a^5}\end{pmatrix}}\in[[g]]$. Therefore taking\;  $a_{21}=\frac{f}{a^2}$; $a_{22}=\frac{b}{a^2}$; $a_{43}=\frac{l}{a^4}$ $a_{44}=\frac{d}{a^4}$ ; $a_{55}=\frac{e}{a^5}$; $a_{33}=\frac{c}{a^3}$ we complete the proof.	
\end{prof}

\subsection{Case of $\mathfrak{g}=A_{5,5}$}
Let's consider $\mathbb{R}^{\times}Aut(\mathfrak{g})$ for this Lie algebra. From the lemma \ref{lem 001}  one can see that $\mathbb{R}\oplus Der(\mathfrak{g})$ which is the Lie algebra of $\mathbb{R}^{\times}Aut(\mathfrak{g})$ contains 
\begin{equation*}
G_6=\left\{{\renewcommand{\arraystretch}{0.4}\begin{pmatrix}
	0&0&0&0&0\\
	a_{21}&0&0&0&0\\
	a_{31}&a_{32}&0&0&0\\
	a_{41}&a_{42}&0&0&0\\
	a_{51}&a_{52}&a_{53}&a_{32}-a_{41}&0
	\end{pmatrix}}\; ;a_{ij}\in \mathbb{R}\right\}. \end{equation*}
Hence, by matrix exponentiation , we obtain

\begin{equation}\label{e009}
F_6=\left\{{\renewcommand{\arraystretch}{0.4}\begin{pmatrix}
	1&0&0&0&0\\
	a_{21}&1&0&0&0\\
	a_{31}&a_{32}&1&0&0\\
	a_{41}&a_{42}&0&1&0\\
	a_{51}&a_{52}&a_{53}&a_{32}-a_{41}&1
	\end{pmatrix}}\; ;a_{ij}\in \mathbb{R} \right\}\subseteq \mathbb{R}^{\times}Aut(\mathfrak{g}). \\\end{equation}
\begin{lemma}\label{lem0061}
	Let $\mathfrak{g}= A_{5,5}$. Then the following $U'_6$ is a set of representatives of $\mathfrak{PM}$	
	\begin{equation}\label{e007}
	U'_6=\left\{{\renewcommand{\arraystretch}{0.4}\begin{pmatrix}
		a&0&0&0&0\\
		0&b&0&0&0\\
		0&0&c&0&0\\
		0&0&l&d&0\\
		0&0&0&f&e
		\end{pmatrix}}\; ; a,b,c,d,e>0,\; l,f\in \mathbb{R}\right\}. \\\end{equation}
\end{lemma}
\begin{prof}
	Take any $g \in GL_5(\mathbb{R})$. Let's find $a,b,c,d,e>0$, $f, l\in \mathbb{R}$ such that: 
	\begin{equation*}
	{\renewcommand{\arraystretch}{0.4}\begin{pmatrix}
		a&0&0&0&0\\
		0&b&0&0&0\\
		0&0&c&0&0\\
		0&0&l&d&0\\
		0&0&0&f&e
		\end{pmatrix}} \in [[g]].
	\end{equation*}
	There exists $k\in O(5)$ such that $gk={\renewcommand{\arraystretch}{0.4}\begin{pmatrix}
		g_{11}&0&0&0&0\\
		g_{21}&g_{22}&0&0&0\\
		g_{31}&g_{32}&g_{33}&0&0\\
		g_{41}&g_{42}&g_{43}&g_{44}&0\\
		g_{51}&g_{52}&g_{53}&g_{54}&g_{55}
		\end{pmatrix}}$	 with $g_{11},g_{22},g_{33},g_{44},g_{55}>0$.\\ Let consider the matrix $\varphi={\renewcommand{\arraystretch}{0.4}\begin{pmatrix}
		1&0&0&0&0\\
		a_{21}&1&0&0&0\\
		a_{31}&a_{32}&1&0&0\\
		a_{41}&a_{42}&0&1&0\\
		a_{51}&a_{52}&a_{53}&a_{32}-a_{41}&1
		\end{pmatrix}}\in F_6\subseteq \mathbb{R}^\times Aut(\mathfrak{g})$ with
	
	$$a_{21}=-\dfrac{g_{21}}{g_{11}}, \quad
	a_{31}=-\frac{g_{31}+a_{32}g_{21}}{g_{11}}, \quad
	a_{32}=-\frac{g_{32}}{g_{22}}, \quad
	a_{41}=-\frac{g_{41}+a_{42}g_{21}}{g_{11}},\quad
	a_{42}=-\frac{g_{42}}{g_{22}},$$
	$$a_{51}=-\frac{g_{51}+(a_{32}-a_{41})g_{41}+a_{53}g_{31}+a_{52}g_{21}}{g_{11}},\quad
	a_{52}=-\frac{g_{52}+a_{53}g_{32}+(a_{32}-a_{41})g_{42}}{g_{22}},$$
	$$a_{53}=-\frac{g_{53}+(a_{32}-a_{41})g_{43}}{g_{33}}.$$
	
Hence $\varphi gk={\renewcommand{\arraystretch}{0.4}\begin{pmatrix}
		g_{11}&0&0&0&0\\
		0&g_{22}&0&0&0\\
		0&0&g_{33}&0&0\\
		0&0&g_{43}&g_{44}&0\\
		0&0&0&(a_{32}-a_{41})g_{44}+g_{54}&g_{55}
		\end{pmatrix}}\in[[g]]$. \\
	One can Take: $a=g_{11},b=g_{22},c=g_{33},d=g_{44}$, $e=g_{55}$, $f=(a_{32}-a_{41})g_{44}+g_{54}$,\; $l=g_{43}$  to conclude the proof. 
\end{prof}

\begin{proposition}\label{propo007}Let $\mathfrak{g}=A_{5,5}$. Then the following $U_6$ is a set of representatives of $\mathfrak{PM}$
	\begin{equation}\label{eq009}	
	U_6=\left\{{\renewcommand{\arraystretch}{0.4}\begin{pmatrix}
		1&0&0&0&0\\
		0&1&0&0&0\\
		0&0&a_{33}&0&0\\
		0&0&a_{43}&a_{44}&0\\
		0&0&0&a_{54}&a_{55}
		\end{pmatrix}}\; ;a_{33},a_{44},a_{55}>0,\; a_{54}\in\mathbb{R}\right\} \\\end{equation}
\end{proposition}

\begin{prof}
	By Lemma \ref{lem0061}, there exists $a,b,c,d,e>0$, $f,l\in \mathbb{R}$ such that $h={\renewcommand{\arraystretch}{0.4}\begin{pmatrix}
		a&0&0&0&0\\
		0&b&0&0&0\\
		0&0&c&0&0\\
		0&0&l&d&0\\
		0&0&0&f&e
		\end{pmatrix}} \in [[g]]$. \\It follows from  \eqref{eq06} that,$A=diag\Big(\frac{1}{a},\frac{1}{b},\frac{1}{b^2},\frac{1}{ab},\frac{1}{ab^2}\Big)
	\in \mathbb{R}^{\times} Aut(\mathfrak{g})$.\\
	Hence, $Ah={\renewcommand{\arraystretch}{0.4}\begin{pmatrix}
		1&0&0&0&0\\
		0&1&0&0&0\\
		0&0&\dfrac{c}{b^2}&0&0\\
		0&0&\frac{l}{ab}&\dfrac{d}{ab}&0\\
		0&0&0&\dfrac{f}{ab^2}&\dfrac{e}{ab^2}\end{pmatrix}}\in[[g]]$. Therefore taking , $a_{54}=\dfrac{f}{ab^2}$; $a_{33}=\dfrac{c}{b^2}$; $a_{43}=\dfrac{l}{ab}$; $a_{44}=\dfrac{d}{ab}$ ; $a_{55}=\dfrac{e}{ab^2}$, we complete the proof.	
\end{prof}

\subsection{Case of $\mathfrak{g}=A_{5,3}$}
Let's consider $\mathbb{R}^{\times}Aut(\mathfrak{g})$ for this Lie algebra. From the lemma \ref{lem 001}  one can see that $\mathbb{R}\oplus Der(\mathfrak{g})$ which is the Lie algebra of $\mathbb{R}^{\times}Aut(\mathfrak{g})$ contains 
\begin{equation*}
G_7=\left\{{\renewcommand{\arraystretch}{0.4}\begin{pmatrix}
	0&0&0&0&0\\
	a_{54}&0&0&0&0\\
	0&a_{32}&0&0&0\\
	a_{41}&a_{42}&a_{32}&0&0\\
	a_{51}&a_{52}&0&a_{54}&0
	\end{pmatrix}}\; ;a_{ij}\in \mathbb{R}\right\}. \end{equation*}
Hence, by matrix exponentiation , we obtain

\begin{equation}\label{e009}
F_7=\left\{{\renewcommand{\arraystretch}{0.4}\begin{pmatrix}
	1&0&0&0&0\\
	a_{54}&1&0&0&0\\
	a_{31}&a_{32}&1&0&0\\
	a_{41}&a_{42}&a_{32}&1&0\\
	a_{51}&a_{52}&a_{53}&a_{54}&1
	\end{pmatrix}}\; ;a_{ij}\in \mathbb{R} \right\}\subseteq \mathbb{R}^{\times}Aut(\mathfrak{g}). \\\end{equation}
\begin{lemma}\label{lem0010}
	Let $\mathfrak{g}= A_{5,3}$. Then the following $U'_7$ is a set of representatives of $\mathfrak{PM}$	
	\begin{equation}\label{e008}
	U'_7=\left\{{\renewcommand{\arraystretch}{0.4}\begin{pmatrix}
		a&0&0&0&0\\
		f&b&0&0&0\\
		0&0&c&0&0\\
		0&0&l&d&0\\
		0&0&0&0&e
		\end{pmatrix}}\; ; a,b,c,d,e>0,\; l,f\in \mathbb{R}\right\}. \\\end{equation}
\end{lemma}
\begin{prof}
	Take any $g \in GL_5(\mathbb{R})$. Let's find  $a,b,c,d,e>0$, $f\in \mathbb{R}$ such that: 
	\begin{equation*}
	{\renewcommand{\arraystretch}{0.4}\begin{pmatrix}
		a&0&0&0&0\\
		f&b&0&0&0\\
		0&0&c&0&0\\
		0&0&l&d&0\\
		0&0&0&0&e
		\end{pmatrix}} \in [[g]].
	\end{equation*}
	There exists $k\in O(5)$ such that $gk={\renewcommand{\arraystretch}{0.4}\begin{pmatrix}
		g_{11}&0&0&0&0\\
		g_{21}&g_{22}&0&0&0\\
		g_{31}&g_{32}&g_{33}&0&0\\
		g_{41}&g_{42}&g_{43}&g_{44}&0\\
		g_{51}&g_{52}&g_{53}&g_{54}&g_{55}
		\end{pmatrix}}$	 with $g_{11},g_{22},g_{33},g_{44},g_{55}>0$.\\ Let consider the matrix $\varphi={\renewcommand{\arraystretch}{0.4}\begin{pmatrix}
		1&0&0&0&0\\
		a_{54}&1&0&0&0\\
		a_{31}&a_{32}&1&0&0\\
		a_{41}&a_{42}&a_{32}&1&0\\
		a_{51}&a_{52}&a_{53}&a_{54}&1
		\end{pmatrix}}\in F_7\subseteq \mathbb{R}^\times Aut(\mathfrak{g})$ with
	$$a_{31}=-\frac{a_{32}g_{21}+g_{31}}{g_{11}}, \quad
	a_{32}=-\frac{g_{32}}{g_{22}}, \quad
	a_{41}=-\frac{g_{41}+a_{32}g_{31}+a_{42}g_{21}}{g_{11}},\quad
	a_{42}=-\frac{g_{42}+a_{32}g_{32}}{g_{22}},$$
	$$a_{51}=-\frac{g_{51}+a_{54}g_{41}+a_{52}g_{54}+a_{53}g_{31}}{g_{11}},\quad
	a_{52}=-\frac{g_{52}+a_{54}g_{42}+a_{53}g_{32}}{g_{22}},\quad
	a_{53}=-\frac{g_{53}+a_{54}g_{43}}{g_{33}},$$
	$$a_{54}=-\frac{g_{54}}{g_{44}}.$$
	
to $\varphi gk={\renewcommand{\arraystretch}{0.4}\begin{pmatrix}
		g_{11}&0&0&0&0\\
		a_{54}g_{11}+g_{21}&g_{22}&0&0&0\\
		0&0&g_{33}&0&0\\
		0&0&a_{32}g_{33}+g_{43}&g_{44}&0\\
		0&0&0&0&g_{55}
		\end{pmatrix}}\in[[g]]$. \\Take $a=g_{11},b=g_{22},c=g_{33},d=g_{44}$, $e=g_{55}$, $l=a_{32}g_{33}+g_{43}$ $a_{54}g_{11}+g_{21}$ and conclude the proof.
\end{prof}

\begin{proposition}\label{propo0010}Let $\mathfrak{g}=A_{5,3}$. Then the following $U_7$ is a set of representatives of $\mathfrak{PM}$
	\begin{equation}\label{eq0011}	
	U_7=\left\{{\renewcommand{\arraystretch}{0.4}\begin{pmatrix}
		1&0&0&0&0\\
		a_{21}&1&0&0&0\\
		0&0&a_{33}&0&0\\
		0&0&a_{43}&a_{44}&0\\
		0&0&0&0&a_{55}
		\end{pmatrix}}\; ;a_{33},a_{44},a_{55}>0,\; a_{43},a_{21}\in\mathbb{R}\right\} \\\end{equation}
\end{proposition}

\begin{prof}
	By Lemma \ref{lem0010}, there exists $a,b,c,d,e>0$, $f,l\in \mathbb{R}$ such that $h={\renewcommand{\arraystretch}{0.4}\begin{pmatrix}
		a&0&0&0&0\\
		f&b&0&0&0\\
		0&0&c&0&0\\
		0&0&l&d&0\\
		0&0&0&0&e
		\end{pmatrix}} \in [[g]]$. \\It follows from  \eqref{eq07} that,$A=diag\Big(\frac{1}{a},\frac{1}{b},\frac{1}{b^2},\frac{1}{ab},\frac{1}{ab^2}\Big)
	\in \mathbb{R}^{\times} Aut(\mathfrak{g})$.\\
	Hence, $Ah={\renewcommand{\arraystretch}{0.4}\begin{pmatrix}
		1&0&0&0&0\\
		\frac{f}{b}&1&0&0&0\\
		0&0&\frac{c}{b^2}&0&0\\
		0&0&\frac{l}{ab}&\frac{d}{ab}&0\\
		0&0&0&0&\frac{e}{ab^2}\end{pmatrix}}\in[[g]]$. Therefore taking, $a_{21}=\frac{f}{b}$; $a_{33}=\frac{c}{b^2}$; $a_{44}=\frac{d}{ab}$ ; $a_{43}=\frac{l}{ab}$  $a_{55}=\frac{e}{ab^2}$ we complete the proof.	
\end{prof}

\subsection{Case of $\mathfrak{g}=A_{5,1}$}
Let's consider $\mathbb{R}^{\times}Aut(\mathfrak{g})$ for this Lie algebra. From the lemma \ref{lem 001}  one can see that $\mathbb{R}\oplus Der(\mathfrak{g})$ which is the Lie algebra of $\mathbb{R}^{\times}Aut(\mathfrak{g})$ contains 
\begin{equation*}
G_8=\left\{{\renewcommand{\arraystretch}{0.4}\begin{pmatrix}
	0&0&0&0&0\\
	a_{21}&0&0&0&0\\
	a_{31}&a_{32}&0&0&0\\
	a_{41}&a_{42}&a_{43}&0&0\\
	a_{51}&a_{52}&a_{53}&a_{32}&0
	\end{pmatrix}}\; ;a_{ij}\in \mathbb{R}\right\}. \end{equation*}
Hence, by matrix exponentiation , we obtain

\begin{equation}\label{eq009}
F_8=\left\{{\renewcommand{\arraystretch}{0.4}\begin{pmatrix}
	1&0&0&0&0\\
	a_{21}&1&0&0&0\\
	a_{31}&a_{32}&1&0&0\\
	a_{41}&a_{42}&a_{41}&1&0\\
	a_{51}&a_{52}&a_{53}&a_{32}&1
	\end{pmatrix}}\; ;a_{ij}\in \mathbb{R} \right\}\subseteq \mathbb{R}^{\times}Aut(\mathfrak{g}). \\\end{equation}
\begin{lemma}\label{lem00010}
	Let $\mathfrak{g}= A_{5,1}$. Then the following $U'_8$ is a set of representatives of $\mathfrak{PM}$	
	\begin{equation}\label{e008}
	U'_8=\left\{{\renewcommand{\arraystretch}{0.4}\begin{pmatrix}
		a&0&0&0&0\\
		0&b&0&0&0\\
		0&0&c&0&0\\
		0&0&0&d&0\\
		0&0&0&f&e
		\end{pmatrix}}\; ; a,b,c,d,e>0,\; f\in \mathbb{R}\right\}. \\\end{equation}
\end{lemma}
\begin{prof}
	Take any $g \in GL_5(\mathbb{R})$. Let's find $a,b,c,d,e>0$, $f\in \mathbb{R}$ such that: 
	\begin{equation*}
	{\renewcommand{\arraystretch}{0.4}\begin{pmatrix}
		a&0&0&0&0\\
		0&b&0&0&0\\
		0&0&c&0&0\\
		0&0&0&d&0\\
		0&0&0&f&e
		\end{pmatrix}} \in [[g]].
	\end{equation*}
	There exists $k\in O(5)$ such that $gk={\renewcommand{\arraystretch}{0.4}\begin{pmatrix}
		g_{11}&0&0&0&0\\
		g_{21}&g_{22}&0&0&0\\
		g_{31}&g_{32}&g_{33}&0&0\\
		g_{41}&g_{42}&g_{43}&g_{44}&0\\
		g_{51}&g_{52}&g_{53}&g_{54}&g_{55}
		\end{pmatrix}}$	 with $g_{11},g_{22},g_{33},g_{44},g_{55}>0$.\\ Let consider the matrix $\varphi={\renewcommand{\arraystretch}{0.4}\begin{pmatrix}
		1&0&0&0&0\\
		a_{21}&1&0&0&0\\
		a_{31}&a_{32}&1&0&0\\
		a_{41}&a_{42}&a_{43}&1&0\\
		a_{51}&a_{52}&a_{53}&a_{32}&1
		\end{pmatrix}} \in F_8\subseteq \mathbb{R}^\times Aut(\mathfrak{g})$ with
	$$a_{21}=-\frac{g_{21}}{g_{11}} \quad
	a_{31}=-\frac{a_{32}g_{21}+g_{31}}{g_{11}} \quad
	a_{32}=-\frac{g_{32}}{g_{22}} \quad
	a_{41}=-\frac{g_{41}+a_{43}g_{31}+a_{42}g_{21}}{g_{11}}\quad a_{42}=-\frac{g_{42}+a_{43}g_{32}}{g_{22}}.$$
	$$a_{43}=-\frac{g_{43}}{g_{33}}\quad
	a_{51}=-\frac{g_{51}+a_{32}g_{41}+a_{52}g_{21}+a_{53}g_{31}}{g_{11}}\quad
	a_{52}=-\frac{g_{52}+a_{32}g_{42}+a_{53}g_{32}}{g_{22}}\quad
	a_{53}=-\frac{g_{53}+a_{32}g_{43}}{g_{33}}.$$		
Hence $\varphi gk={\renewcommand{\arraystretch}{0.4}\begin{pmatrix}
		g_{11}&0&0&0&0\\
		0&g_{22}&0&0&0\\
		0&0&g_{33}&0&0\\
		0&0&0&g_{44}&0\\
		0&0&0&a_{32}g_{44}+g_{54}&g_{55}
		\end{pmatrix}}\in[[g]]$. \\ One can Take $a=g_{11},b=g_{22},c=g_{33},d=g_{44}$, $e=g_{55}$, $f=a_{32}g_{44}+g_{54}$ to conclude the proof.
\end{prof}

\begin{proposition}\label{propo00010}Let $\mathfrak{g}=A_{5,1}$. Then the following $U_8$ is a set of representatives of $\mathfrak{PM}$
	\begin{equation}\label{eq00011}	
	U_8=\left\{{\renewcommand{\arraystretch}{0.4}\begin{pmatrix}
		1&0&0&0&0\\
		0&1&0&0&0\\
		0&0&1&0&0\\
		0&0&0&a_{44}&0\\
		0&0&0&a_{54}&a_{55}
		\end{pmatrix}}\; ;a_{44},a_{55}>0,\; a_{54}\in\mathbb{R}\right\}. \\\end{equation}
\end{proposition}

\begin{prof}
	By Lemma \ref{lem00010}, there exists $a,b,c,d,e>0$, $f\in \mathbb{R}$ such that $h={\renewcommand{\arraystretch}{0.4}\begin{pmatrix}
		a&0&0&0&0\\
		0&b&0&0&0\\
		0&0&c&0&0\\
		0&0&0&d&0\\
		0&0&0&f&e
		\end{pmatrix}} \in [[g]]$. \\It follows from  \eqref{eq08} that,$A=diag\Big(\frac{1}{a},\frac{1}{b},\frac{1}{c},\frac{1}{ab},\frac{1}{ac}\Big)
	\in \mathbb{R}^{\times} Aut(\mathfrak{g})$.\\
	Hence,$Ah={\renewcommand{\arraystretch}{0.4}\begin{pmatrix}
		1&0&0&0&0\\
		0&1&0&0&0\\
		0&0&1&0&0\\
		0&0&0&\frac{d}{ab}&0\\
		0&0&0&\frac{f}{ac}&\frac{e}{ac}\end{pmatrix}}\in[[g]]$. Therefore taking ,$a_{54}=\frac{f}{ac}$; $a_{44}=\frac{d}{ab}$ ; $a_{55}=\frac{e}{ac}$ we complete the proof.	
\end{prof}

\subsection{Case of $\mathfrak{g}=A_{5,2}$}
Let's consider $\mathbb{R}^{\times}Aut(\mathfrak{g})$ for this Lie algebra. From the lemma \ref{lem 001}  one can see that $\mathbb{R}\oplus Der(\mathfrak{g})$ which is the Lie algebra of $\mathbb{R}^{\times}Aut(\mathfrak{g})$ contains 
\begin{equation*}
G_9=\left\{{\renewcommand{\arraystretch}{0.4}\begin{pmatrix}
	0&0&0&0&0\\
	a_{21}&0&0&0&0\\
	a_{31}&a_{43}&0&0&0\\
	a_{41}&a_{42}&a_{43}&0&0\\
	a_{51}&a_{52}&a_{42}&a_{54}&0
	\end{pmatrix}}\; ;a_{ij}\in \mathbb{R}\right\}. \end{equation*}
Hence, by matrix exponentiation, we obtaints

\begin{equation}\label{e009}
F_9=\left\{{\renewcommand{\arraystretch}{0.4}\begin{pmatrix}
	1&0&0&0&0\\
	a_{21}&1&0&0&0\\
	a_{31}&a_{43}&1&0&0\\
	a_{41}&a_{42}&a_{43}&1&0\\
	a_{51}&a_{52}&a_{53}&a_{54}&1
	\end{pmatrix}}\; ;a_{ij}\in \mathbb{R} \right\}\subseteq \mathbb{R}^{\times}Aut(\mathfrak{g}). \\\end{equation}
\begin{lemma}\label{lem009}
	Let $\mathfrak{g}= A_{5,2}$. Then the following $U'_9$ is a set of representatives of $\mathfrak{PM}$	
	\begin{equation}\label{e008}
	U'_9=\left\{{\renewcommand{\arraystretch}{0.4}\begin{pmatrix}
		a&0&0&0&0\\
		0&b&0&0&0\\
		0&f&c&0&0\\
		0&0&0&d&0\\
		0&0&0&0&e
		\end{pmatrix}}\; ; a,b,c,d,e>0,\; f\in \mathbb{R}\right\}. \\\end{equation}
\end{lemma}
\begin{prof}
	Take any $g \in GL_5(\mathbb{R})$. Let's fine  $a,b,c,d,e>0$, $f\in \mathbb{R}$ such that: 
	\begin{equation*}
	{\renewcommand{\arraystretch}{0.4}\begin{pmatrix}
		a&0&0&0&0\\
		0&b&0&0&0\\
		0&f&c&0&0\\
		0&0&0&d&0\\
		0&0&0&0&e		\end{pmatrix}} \in [[g]].
	\end{equation*}
	There exists $k\in O(5)$ such that $gk={\renewcommand{\arraystretch}{0.4}\begin{pmatrix}
		g_{11}&0&0&0&0\\
		g_{21}&g_{22}&0&0&0\\
		g_{31}&g_{32}&g_{33}&0&0\\
		g_{41}&g_{42}&g_{43}&g_{44}&0\\
		g_{51}&g_{52}&g_{53}&g_{54}&g_{55}
		\end{pmatrix}}$	 with $g_{11},g_{22},g_{33},g_{44},g_{55}>0$.\\ Let consider the matrix	
	$\varphi={\renewcommand{\arraystretch}{0.4}\begin{pmatrix}
		1&0&0&0&0\\
		a_{21}&1&0&0&0\\
		a_{31}&a_{43}&1&0&0\\
		a_{41}&a_{42}&a_{43}&1&0\\
		a_{51}&a_{52}&a_{53}&a_{54}&1
		\end{pmatrix}} \in F_9\subseteq \mathbb{R}^\times Aut(\mathfrak{g})$ with
	$$a_{21}=-\frac{g_{21}}{g_{11}}, \quad
	a_{31}=-\frac{a_{43}g_{21}+g_{31}}{g_{11}}, \quad
	a_{41}=-\frac{g_{41}+a_{42}g_{21}+a_{43}g_{31}}{g_{11}},\quad
	a_{42}=-\frac{g_{42}+a_{43}g_{32}}{g_{22}},\quad
	a_{43}=-\frac{g_{43}}{g_{33}},$$
	$$a_{51}=-\frac{g_{51}+a_{54}g_{41}+a_{52}g_{21}+a_{53}g_{31}}{g_{11}},\quad
	a_{52}=-\frac{g_{52}+a_{54}g_{42}+a_{53}g_{32}}{g_{22}},\quad
	a_{53}=-\frac{g_{53}+a_{54}g_{43}}{g_{44}},\quad
	a_{54}=-\frac{g_{54}}{g_{44}}$$
		
	Hence $\varphi gk={\renewcommand{\arraystretch}{0.4}\begin{pmatrix}
		g_{11}&0&0&0&0\\
		0&g_{22}&0&0&0\\
		0&a_{43}g_{22}+g_{32}&g_{33}&0&0\\
		0&0&0&g_{44}&0\\
		0&0&0&0&g_{55}
		\end{pmatrix}}\in[[g]]$. \\One can Take $a=g_{11},b=g_{22},c=g_{33},d=g_{44}$, $e=g_{55}$, $f=a_{43}g_{22}+g_{32}$ to conclude the proof.
\end{prof}

\begin{proposition}\label{propo009}Let's $\mathfrak{g}=A_{5,2}$. Then the following $U_9$ is a set of representatives of $\mathfrak{PM}$
	\begin{equation}\label{eq009}	
	U_9=\left\{{\renewcommand{\arraystretch}{0.4}\begin{pmatrix}
		1&0&0&0&0\\
		0&1&0&0&0\\
		0&a_{32}&a_{33}&0&0\\
		0&0&0&a_{44}&0\\
		0&0&0&0&a_{55}
		\end{pmatrix}}\; ;a_{33},a_{44},a_{55}>0,\;a_{43}\in\mathbb{R}\right\}. \\\end{equation}
\end{proposition}

\begin{prof}
	By Lemma \ref{lem009}, there exists $a,b,c,d,e>0$, $f,l\in \mathbb{R}$ such that $h={\renewcommand{\arraystretch}{0.4}\begin{pmatrix}
		a&0&0&0&0\\
		0&b&0&0&0\\
		0&f&c&0&0\\
		0&0&0&d&0\\
		0&0&0&0&e
		\end{pmatrix}} \in [[g]]$. \\It follows from  \eqref{eq09} that,$A=diag\Big(\frac{1}{a},\frac{1}{b},\frac{1}{ab},\frac{1}{a^2b},\frac{1}{ab^2}\Big)
	\in \mathbb{R}^{\times} Aut(\mathfrak{g})$.\\
	Hence,$Ah={\renewcommand{\arraystretch}{0.4}\begin{pmatrix}
		1&0&0&0&0\\
		0&1&0&0&0\\
		0&\frac{f}{ab}&\frac{c}{ab}&0&0\\
		0&0&0&\frac{d}{a^2b}&0\\
		0&0&0&0&\frac{e}{ab^2}\end{pmatrix}}\in[[g]]$. Therefore taking, $a_{32}=\frac{f}{ab}$; $a_{33}=\frac{c}{ab}$; $a_{44}=\frac{d}{a^2b}$ ; $a_{55}=\frac{e}{ab^2}$ we complete the proof.	
\end{prof}

\section{Milnor-type theorems for left-invariant Riemannian metrics on $5$-dimensional nilpotent Lie groups}

Throughout this section we apply the Milnor-type theorems on the nilpotent Lie algebras cited in the previous table  to describe inner products.
\subsection{Case of $\mathfrak{g}=5A_{1}$}
\begin{proposition}\label{propoo}
	For every inner product $\langle,\rangle$ on $\mathfrak{g}=5A_{1}$, there exist $\eta> 0$,
	and an orthonormal basis $\mathcal{B}_1=\{v_1, v_2, v_3,v_4,v_5\}$ with respect to $\eta\langle,\rangle$, such that all the commutators are zero
\end{proposition}
\begin{prof}
The theoreme\ref{theo01} gives $k>0$ and $\varphi\in Aut(\mathfrak{g})$ such that the basis $\mathcal{B}_1=\{\varphi e_1,\varphi e_2,\varphi e_3,\varphi e_4,\varphi e_5\}$ is orthormal whith respect to $k\langle,\rangle$. We can take $\eta=k$ to conclude the proof.
\end{prof}
\subsection{Case of $\mathfrak{g}=A_{5,4}$}
\begin{proposition}\label{prop00001}
	For every inner product $\langle,\rangle$ on $\mathfrak{g}=A_{5,4}$, there exist $\alpha\in\mathbb{R}$, $\beta, \gamma>0$, $\eta> 0$,
	and an orthonormal basis $\mathcal{B}_2=\{v_1, v_2, v_3,v_4,v_5\}$ with respect to $\eta\langle,\rangle$, such that non zero commutators are given by
	\begin{equation}
	[v_1,v_3]=\alpha v_5;\quad [v_1,v_4]=\beta v_5;\quad [v_2,v_3]=\gamma v_5 \end{equation}
\end{proposition}
\begin{prof}
	Let $\langle,\rangle$ be an inner product on $\mathfrak{g}=A_{5,4}$. 
	
The theorem \ref{theo01} gives $h  ={\renewcommand{\arraystretch}{0.4}\begin{pmatrix}
		1&0&0&0&0\\
		a_{21}&1&0&0&0\\
		0&0&1&0&0\\
		0&0&0&a_{44}&0\\
		0&0&0&0&a_{55}
		\end{pmatrix}} 
	\in U_2$ where $U_2$ is given by proposition \ref{propo00001}, $\varphi\in Aut(\mathfrak{g})$ and $k>0$ such that $\mathcal{B}_2=\{\varphi he_1,\cdots,\varphi he_5\}$ is an orthonormal basis of $\mathfrak{g}$ whith respect to $k\langle,\rangle$.Then we have:
	\begin{equation*}
	[v_1,v_4]=a_{44}\varphi[he_1,he_4]=a_{44}\varphi[e_1,e_4]=\frac{a_{44}}{a_{55}}\varphi he_{5}=\frac{a_{44}}{a_{55}}v_5;
	\end{equation*} 
	\begin{equation*}
	[v_1,v_3]=\varphi[he_1,he_3]=a_{21}\varphi[e_2,e_3]=\frac{a_{21}}{a_{55}}\varphi he_{5}=\frac{a_{21}}{a_{55}}v_5;
	\end{equation*}
	\begin{equation*}
	[v_2,v_3]=\varphi[he_2,he_3]=\varphi[e_2,e_3]=\frac{1}{a_{55}}\varphi he_{5}=\frac{1}{a_{55}}v_5.
	\end{equation*}
	with $\alpha=\frac{a_{21}}{a_{55}}$,\;$\beta=\frac{a_{44}}{a_{55}}$,\;$\gamma=\frac{1}{a_{55}}$,\;$\eta=k$ we conclude de proof.
\end{prof}

\subsection{Case of $\mathfrak{g}=A_{3,1}\oplus 2A_{1}$}
\begin{proposition}\label{prop00002}
	For every inner product $\langle,\rangle$ on $\mathfrak{g}=A_{3,1}\oplus 2A_{1}$, there exist $\alpha>0$, $\eta> 0$,
	and an orthonormal basis $\mathcal{B}_3=\{v_1, v_2, v_3,v_4,v_5\}$ with respect to $\eta\langle,\rangle$, such that non zero commutators are given by
	\begin{equation}
	[v_1,v_2]=\alpha v_5 \end{equation}
\end{proposition}
\begin{prof}
	Let $\langle,\rangle$ be an inner product on $\mathfrak{g}=A_{3,1}\oplus 2A_{1}$.
The theorem \ref{theo01} gives $h  =diag(1,1,1,1,a_{55})
	\in U_3$, where $U_3$ is given by proposition \ref{propo002}, $\varphi\in Aut(\mathfrak{g})$ and $k>0$ such that $\mathcal{B}_3=\{\varphi he_1,\cdots,\varphi he_5\}$ is an orthonormal basis of $\mathfrak{g}$ whith respect to $k\langle,\rangle$. Then we have: 
	\begin{equation*}
	[v_1,v_2]=\varphi[he_1,he_2]=\varphi[e_1,e_2]=\frac{1}{a_{55}}\varphi he_{5}=\frac{1}{a_{55}}v_5.
	\end{equation*}
	With $\alpha=a_{55}$,\;$\eta=k$ we conclude de proof.
\end{prof}

\subsection{Case of $\mathfrak{g}=A_{4,1}\oplus A_{1}$}
\begin{proposition}\label{prpo6}
	For every inner product $\langle,\rangle$ on $\mathfrak{g}=A_{4,1}\oplus A_{1}$, there exist $\alpha,\beta> 0$, $\gamma\in\mathbb{R}$, $\eta> 0$,
	and an orthonormal basis $\mathcal{B}_4=\{v_1, v_2, v_3,v_4,v_5\}$ with respect to $\eta\langle,\rangle$, such that non zero commutators are given by
	\begin{equation}
	[v_1,v_2]=\alpha v_3 +\gamma v_5;\quad [v_1,v_3]=\beta v_5\end{equation}
	Or
	\begin{equation}
	[v_1,v_2]=\alpha v_3 +\gamma v_4;\quad [v_1,v_3]=\beta v_5.\end{equation}	
\end{proposition}
\begin{prof}
	Let $\langle,\rangle$ be an inner product on $\mathfrak{g}=A_{4,1}\oplus A_{1}$.
\\
The theorem \ref{theo01} gives $h ={\renewcommand{\arraystretch}{0.4}\begin{pmatrix}
		1&0&0&0&0\\
		0&1&0&0&0\\
		0&0&a&0&0\\
		0&0&0&1&0\\
		0&0&c&0&b
		\end{pmatrix}}\in U^1_4$ or $h ={\renewcommand{\arraystretch}{0.4}\begin{pmatrix}
		1&0&0&0&0\\
		0&1&0&0&0\\
		0&0&a&0&0\\
		0&0&c&1&0\\
		0&0&0&0&b
		\end{pmatrix}}\in U^2_4$ where $U^1_4$ and $U^2_4$ are given by proposition \ref{propo003}, $\varphi\in Aut(\mathfrak{g})$ and $k>0$ such that $\mathcal{B}_4=\{\varphi he_1,\cdots,\varphi he_5\}$ is an orthonormal basis of $\mathfrak{g}$ whith respect to $k\langle,\rangle$. Then we have:
	\begin{enumerate}
		\item [1)]\underline{First case:} if $h\in U^1_4$
		\begin{equation*}
		[v_1,v_2]=\varphi[he_1,he_2]=\varphi[e_1,e_2]=\varphi e_3=\frac{1}{a}v_3-\frac{c}{ab}v_5;
		\end{equation*}
		\begin{equation*}
		[v_1,v_3]=\varphi[he_1,he_3]=\varphi[e_1,ae_3]=a\varphi e_5=\frac{a}{b}v_5.
		\end{equation*}
		With $\alpha=\frac{1}{a}$,\;$\beta=\frac{a}{b}$,\; $\gamma=-\frac{c}{ab}$,\;$\eta=k$ we conclude the proof.
		\item [2)] \underline{Second case:}
		if $h\in U^2_4$
		\begin{equation*}
		[v_1,v_2]=\varphi[he_1,he_2]=\varphi[e_1,e_2]=\varphi e_3=\frac{1}{a}v_3-\frac{c}{a}v_4;
		\end{equation*}
		\begin{equation*}
		[v_1,v_3]=\varphi[he_1,he_3]=\varphi[e_1,ae_3]=a\varphi e_5=\frac{a}{b}v_5.
		\end{equation*}
	\end{enumerate}
	With $\alpha=\frac{1}{a}$,\;$\beta=\frac{a}{b}$,\; $\gamma=-\frac{c}{a}$,\;$\eta=k$ we conclude the proof.
\end{prof}

\subsection{Case of $\mathfrak{g}=A_{5,6}$}
\begin{proposition}\label{prpo7}
	For every inner product $\langle,\rangle$ on $\mathfrak{g}=A_{5,6}$, there exist $\alpha<0$ ,$\gamma,\varepsilon,\sigma> 0$;\; $\beta,\delta\in\mathbb{R}$, $\eta> 0$,
	and an orthonormal basis $\mathcal{B}_5=\{v_1, v_2, v_3,v_4,v_5\}$ with respect to $\eta\langle,\rangle$, such that non zero commutators are given by
	\begin{equation}
	[v_1,v_2]=\alpha v_3+\beta v_4;\quad [v_1,v_3]=\gamma v_4+\delta v_5;\quad [v_1,v_4]=\varepsilon v_5,\quad [v_2,v_3]=\sigma v_5. \end{equation}
\end{proposition}
\begin{prof}
	Let $\langle,\rangle$ be an inner product on $\mathfrak{g}=A_{5,6}$.
\\
	The theorem \ref{theo01} gives $h =\begin{pmatrix}
	1&0&0&0&0\\
	a_{21}&a_{22}&0&0&0\\
	0&0&a_{33}&0&0\\
	0&0&a_{43}&a_{44}&0\\
	0&0&0&0&a_{55}
	\end{pmatrix}\in U_5$ where $U_5$ is given by proposition \ref{propo004}, $\varphi\in Aut(\mathfrak{g})$ and $k>0$ such that $\mathcal{B}_5=\{\varphi he_1,\cdots,\varphi he_5\}$ is an orthonormal basis of $\mathfrak{g}$ whith respect to $k\langle,\rangle$. Then we have:
	\begin{equation*}
	[v_1,v_2]=\varphi[he_1,he_2]=-\frac{a_{22}}{a_{33}}v_3+\frac{a_{22}a_{43}}{a_{33}a_{44}}v_4;\quad[v_1,v_3]=\varphi[he_1,he_3]=\frac{a_{33}}{a_{44}}v_4+\Big(\frac{a_{44}+a_{21}a_{33}}{a_{55}}\Big)v_5; 
	\end{equation*}
	\begin{equation*}
	[v_1,v_4]=\varphi[he_1,he_4]=\frac{a_{44}}{a_{55}}v_5;\quad[v_2,v_3]=\varphi[he_2,he_3]=\frac{a_{22}a_{33}}{a_{55}}v_5.
	\end{equation*}
	With $\alpha=-\frac{a_{22}}{a_{33}}$,\; $\beta=\frac{a_{22}a_{43}}{a_{33}a_{44}}$ ,\;$\gamma=\frac{a_{33}}{a_{44}}$,\; $\delta=\frac{a_{44}+a_{21}a_{33}}{a_{55}}$,\;  $\varepsilon=\frac{a_{44}}{a_{55}}$,\;$\sigma=\frac{a_{22}a_{33}}{a_{55}}$,\;$\eta=k$ we conclude de proof.
\end{prof}

\subsection{Case of $\mathfrak{g}=A_{5,5}$}
\begin{proposition}\label{prpo8}
	For every inner product $\langle,\rangle$ on $\mathfrak{g}=A_{5,5}$, there exist $\alpha,\gamma,\varepsilon> 0$;\; $\beta,\delta\in\mathbb{R}$, $\eta> 0$
	and an orthonormal basis $\mathcal{B}_6=\{v_1, v_2, v_3,v_4,v_5\}$ with respect to $\eta\langle,\rangle$, such that non zero commutators are given by
	\begin{equation}
	[v_1,v_2]=\alpha v_4+\beta v_5;\quad [v_1,v_3]=\gamma v_5;\quad [v_2,v_3]=\delta v_5,\quad [v_2,v_4]=\varepsilon v_5. \end{equation}
\end{proposition}
\begin{prof}
	Let $\langle,\rangle$ be an inner product on $\mathfrak{g}=A_{5,5}$.
\\
The theorem \ref{theo01} gives $h =\begin{pmatrix}
	1&0&0&0&0\\
	0&1&0&0&0\\
	0&0&a_{33}&0&0\\
	0&0&a_{43}&a_{44}&0\\
	0&0&0&a_{54}&a_{55}
	\end{pmatrix}\in U_6$  where $U_6$ is given by proposition \ref{propo007},\\ $\varphi\in Aut(\mathfrak{g})$ and $k>0$ such that $\mathcal{B}_6=\{\varphi he_1,\cdots,\varphi he_5\}$ is an orthonormal basis of $\mathfrak{g}$ whith respect to $k\langle,\rangle$. Then we have:
	\begin{equation*}
	[v_1,v_2]=\varphi[he_1,he_2]=\frac{1}{a_{44}}v_4-\frac{a_{54}}{a_{44}a_{55}}v_5;\; [v_1,v_3]=\varphi[he_1,he_3]=\frac{a_{33}}{a_{55}}v_5;
	\end{equation*} 
	\begin{equation*}
	[v_2,v_3]=\varphi[he_2,he_3]=\frac{a_{43}}{a_{55}}v_5; \;[v_2,v_4]=\varphi[he_2,he_4]=\frac{a_{44}}{a_{55}}v_5.\end{equation*}
	With $\alpha=\frac{1}{a_{44}}$,\;$\beta=-\frac{a_{54}}{a_{44}a_{55}}$,\; $\gamma=\frac{a_{33}}{a_{55}}$,\; $\delta=\frac{a_{43}}{a_{55}}$,\; $\varepsilon=\frac{a_{44}}{a_{55}}$,\; $\eta=k$ we conclude de proof.
\end{prof}

\subsection{Case of $\mathfrak{g}=A_{5,3}$}
\begin{proposition}\label{prpo08}
	For every inner product $\langle,\rangle$ on $\mathfrak{g}=A_{5,3}$, there exist $\alpha,\gamma,\varepsilon> 0$;\; $\beta,\delta\in\mathbb{R}$, $\eta> 0$,
	and an orthonormal basis $\mathcal{B}_7=\{v_1, v_2, v_3,v_4,v_5\}$ with respect to $\eta\langle,\rangle$, such that non zero commutators are given by
	\begin{equation}
	[v_1,v_2]=\alpha v_3+\beta v_4;\quad [v_1,v_3]=\gamma v_4+\delta v_5;\quad [v_2,v_3]=\varepsilon v_5. \end{equation}
\end{proposition}
\begin{prof}
	Let $\langle,\rangle$ be an inner product on $\mathfrak{g}=A_{5,3}$.
\\
The theorem \ref{theo01} gives $h =\begin{pmatrix}
	1&0&0&0&0\\
	a_{21}&1&0&0&0\\
	0&0&a_{33}&0&0\\
	0&0&a_{43}&a_{44}&0\\
	0&0&0&0&a_{55}
	\end{pmatrix}\in U_7$  where $U_7$ is given by proposition \ref{propo0010}, $\varphi\in Aut(\mathfrak{g})$ and $k>0$ such that $\mathcal{B}_7=\{\varphi he_1,\cdots,\varphi he_5\}$ is an orthonormal basis of $\mathfrak{g}$ whith respect to $k\langle,\rangle$. Then we have:
	\begin{equation*}
	[v_1,v_2]=\varphi[he_1,he_2]=\frac{1}{a_{33}}v_3-\frac{a_{43}}{a_{33}a_{44}}v_4;\;[v_1,v_3]=\varphi[he_1,he_3]=\frac{a_{33}}{a_{44}}v_4+\frac{a_{21}a_{33}}{a_{55}}v_5.
	\end{equation*} 
	\begin{equation*}
	[v_2,v_3]=\varphi[he_2,he_3]=\frac{a_{33}}{a_{55}}v_5,\; 
	\end{equation*}
	With $\alpha=\frac{1}{a_{33}}$,\;$\beta=-\frac{a_{43}}{a_{33}a_{44}}$,\; $\gamma=\frac{a_{33}}{a_{44}}$,\; $\delta=\frac{a_{21}a_{33}}{a_{55}}$,\; $\varepsilon=\frac{a_{33}}{a_{55}}$,\; $\eta=k$ we conclude de proof.
\end{prof}

\subsection{Case of $\mathfrak{g}=A_{5,1}$}
\begin{proposition}\label{prpo008}
	For every inner product $\langle,\rangle$ on $\mathfrak{g}=A_{5,1}$, there exist $\alpha,\gamma> 0$;\; $\beta\in\mathbb{R}$, $\eta> 0$,
	and an orthonormal basis $\mathcal{B}_8=\{v_1, v_2, v_3,v_4,v_5\}$ with respect to $\eta\langle,\rangle$, such that non zero commutators are given by
	\begin{equation}
	[v_1,v_2]=\alpha v_4+\beta v_5;\quad [v_1,v_3]=\gamma v_4. \end{equation}
\end{proposition}
\begin{prof}
	Let $\langle,\rangle$ be an inner product on $\mathfrak{g}=A_{5,1}$.
\\
	The theorem \ref{theo01} gives of $h =\begin{pmatrix}
	1&0&0&0&0\\
	0&1&0&0&0\\
	0&0&1&0&0\\
	0&0&0&a_{44}&0\\
	0&0&0&a_{54}&a_{55}
	\end{pmatrix}\in U_8$  where $U_8$ is given by proposition \ref{propo00010}, $\varphi\in Aut(\mathfrak{g})$ and $k>0$ such that $\mathcal{B}_8=\{\varphi he_1,\cdots,\varphi he_5\}$ is an orthonormal basis of $\mathfrak{g}$ whith respect to $k\langle,\rangle$. Then we have:
	\begin{equation*}
	[v_1,v_2]=\varphi[he_1,he_2]=\frac{1}{a_{44}}v_4-\frac{a_{54}}{a_{44}a_{55}}v_5;\;	[v_1,v_3]=\varphi[he_1,he_3]=\frac{1}{a_{55}}v_5.
	\end{equation*} 
	With $\alpha=\frac{1}{a_{44}}$,\;$\beta=-\frac{a_{54}}{a_{44}a_{55}}$,\; $\gamma=\frac{1}{a_{55}}$,\; $\eta=k$ we conclude de proof.
\end{prof}

\subsection{Case of $\mathfrak{g}=A_{5,2}$}
\begin{proposition}\label{prpo10}
	For every inner product $\langle,\rangle$ on $\mathfrak{g}=A_{5,2}$, there exist $\alpha,\gamma,\delta> 0$;\; $\beta\in\mathbb{R}$, $\eta> 0$,
	and an orthonormal basis $\mathcal{B}_9=\{v_1, v_2, v_3,v_4,v_5\}$ with respect to $\eta\langle,\rangle$, such that non zero commutators are given by
	\begin{equation}
	[v_1,v_2]=\alpha v_3+\beta v_4;\quad [v_1,v_3]=\gamma v_4;\quad [v_1,v_4]=\delta v_5. \end{equation}
\end{proposition}
\begin{prof}
	Let $\langle,\rangle$ be an inner product on $\mathfrak{g}=A_{5,2}$.
%
\\
The theorem \ref{theo01} gives $h =\begin{pmatrix}
	1&0&0&0&0\\
	0&1&0&0&0\\
	0&a_{32}&a_{33}&0&0\\
	0&0&0&a_{44}&0\\
	0&0&0&0&a_{55}
	\end{pmatrix}\in U_9$  where $U_9$ is given by proposition \ref{propo009}, $\varphi\in Aut(\mathfrak{g})$ and $k>0$ such that $\mathcal{B}_9=\{\varphi he_1,\cdots,\varphi he_5\}$ is an orthonormal basis of $\mathfrak{g}$ whith respect to $k\langle,\rangle$. Then we have:
	\begin{equation*}
	[v_1,v_2]=\varphi[he_1,he_2]=\frac{1}{a_{33}}v_3+\frac{a_{32}}{a_{44}}v_4; \;[v_1,v_3]=\varphi[he_1,he_3]=\frac{a_{33}}{a_{44}}v_4;\;[v_1,v_4]=\varphi[he_1,he_4]=\frac{a_{44}}{a_{55}}v_5.\;
	\end{equation*} 
	With $\alpha=\frac{1}{a_{33}}$,\;$\beta=\frac{a_{32}}{a_{44}}$,\; $\gamma=\frac{a_{33}}{a_{44}}$,\; $\delta=\frac{a_{44}}{a_{55}}$,\;$\eta=k$ we conclude de proof.
\end{prof}
\section{Ricci curvature equation on 5-dimensional nilpotent Lie groups.}

\subsection{Ricci curvature of left-invariant metrics on Lie groups.}
In this subsection, we recall some tools on Ricci curvature of left invariant metrics on Lie groups see\cite{buc} for more details.
Let $(G,g)$ be a Riemaniann Lie group with Lie algebra $(\mathfrak{g},[.,.])$.
For any $u\in\mathfrak{g}$, we denote by $ad_u, \mathsf{J}_u:\mathfrak{g}\longrightarrow\mathfrak{g}$ the endomorphisms given by $ad_uv=[u,v]$ and, $\mathsf{J}_uv=ad^*_vu$ where $ad^*_v$ is the transpose of $ad_v$. Let $B : \mathfrak{g}\times \mathfrak{g}\longrightarrow\mathfrak{g}$ the Killing form given by
\begin{equation*}
B(u,v)=tr(ad_u\circ ad_v) \quad \forall u,v\in\mathfrak{g}.
\end{equation*}
The mean curvature vector $H$ on $\mathfrak{g}$ is the vector given by
\begin{equation*}
\langle H,u\rangle = tr(adu), \forall u \in \mathfrak{g} \quad where \langle.,.\rangle=g(e).
\end{equation*}

\begin{theorem}\cite{buc}
	For any $u,v\in\mathfrak{g}$, 
	\begin{equation}
	ric(u,v)=-\frac{1}{2}B(u,v)-\frac{1}{2}tr(ad_u\circ ad^*_v)-\frac{1}{4}tr(\mathsf{J}_u\circ\mathsf{J}_v)-\frac{1}{2}(\langle ad_Hu, v\rangle+\langle ad_Hv, u\rangle).
	\end{equation}
\end{theorem}
\begin{corollary}\cite{buc}
	Let $(G, g)$ be a Riemannian nilpotent Lie group. Then its Ricci
	curvature at $e$ is given by
	\begin{equation}\label{riccicur}
	ric(u,v)=-\frac{1}{2}tr(ad_u\circ ad^*_v)-\frac{1}{4}tr(\mathsf{J}_u\circ\mathsf{J}_v).
	\end{equation}
\end{corollary}
\begin{remark}
	If $u\in Z(\mathfrak{g})$, $ad_u=0$ and when $u\in\mathfrak{D}^{\perp}(\mathfrak{g})$, $\mathsf{J}_u=0$. 
\end{remark}
\subsection{Ricci curvature of left-invriant metrics on 5-dimensional nilpotent Lie groups.}
In this subsection, we compute the Ricci curvature of left-invariant metrics on $5$-dimensionnal nilpotent Lie groups.
By identifying the endomorphisms $ad_u$ and $\mathsf{J}_u$ with their matrices in the basis obtained in section $4$.

\subsubsection{Case of $\mathfrak{g}=5A_{1}$}
In the basis $\mathcal{B}_1=\{v_1,v_2,v_3,v_4,v_5\}$ of proposition \ref{propoo}, a direct computation using \eqref{riccicur} gives

$$[Ric_{\eta\langle\cdotp,\cdotp\rangle}]=	[Ric_{\langle\cdotp,\cdotp\rangle}]=0$$

\subsubsection{Case of $\mathfrak{g}=A_{5,4}$}
In the basis $\mathcal{B}_2=\{v_1,v_2,v_3,v_4,v_5\}$ in proposition \ref{prop00001}, a direct computation using \eqref{riccicur} gives
\begin{equation}
[Ric_{\eta\langle\cdotp,\cdotp\rangle}]=	[Ric_{\langle\cdotp,\cdotp\rangle}]=-\frac{1}{2}{\renewcommand{\arraystretch}{0.4}\begin{pmatrix}
	\alpha^2+\beta^2&\alpha\gamma&0&0&0\\
	\alpha\gamma&\gamma^2&0&0&0\\
	0&0&\alpha^2+\gamma^2&\alpha\beta&0\\
	0&0&\alpha\beta&\beta^2&0\\
	0&0&0&0&-\alpha^2-\beta^2-\gamma^2
	\end{pmatrix}}. \label{eqn00005}
\end{equation}

\subsubsection{Case of $\mathfrak{g}=A_{3,1}\oplus 2A_1$}
In the basis $\mathcal{B}_3=\{v_1,v_2,v_3,v_4,v_5\}$ in proposition \ref{prop00002}, a direct computation using \eqref{riccicur} gives
\begin{equation}
[Ric_{\eta\langle\cdotp,\cdotp\rangle}]=	[Ric_{\langle\cdotp,\cdotp\rangle}]=-\frac{1}{2}(\alpha^2,\alpha^2,0,0,-\alpha^2) \label{eqan000005}
\end{equation}

\subsubsection{Case of $\mathfrak{g}=A_{4,1}\oplus A_1$}
In the basis $\mathcal{B}_4=\{v_1,v_2,v_3,v_4,v_5\}$in proposition \ref{prpo6}, a direct computation using \eqref{riccicur} gives
\begin{itemize}
	\item[1.]\underline{First case}
 we have de following expression of the Ricci curvature in the basis $\mathcal{B}_4=\{v_1,v_2,v_3,v_4,v_5\}$
	\begin{equation}
	[Ric_{\eta\langle\cdotp,\cdotp\rangle}]=[Ric_{\langle\cdotp,\cdotp\rangle}]=-\frac{1}{2}{\renewcommand{\arraystretch}{0.4}\begin{pmatrix}
		\alpha^2+\beta^2+\gamma^2&0&0&0&0\\
		0&\alpha^2+\gamma^2&\beta\gamma&0&0\\
		0&\beta\gamma&\beta^2-\alpha^2&0&-\alpha\gamma\\
		0&0&0&0&0\\
		0&0&-\alpha\gamma&0&-\beta^2-\gamma^2
		\end{pmatrix}}. \label{eqn4}
	\end{equation}
	\item[2.]\underline{Second case}
	\begin{equation}
	[Ric_{\eta\langle\cdotp,\cdotp\rangle}]= [Ric_{\langle\cdotp,\cdotp\rangle}]=-\frac{1}{2}{\renewcommand{\arraystretch}{0.4}\begin{pmatrix}
		\alpha^2+\beta^2+\gamma^2&0&0&0&0\\
		0&\alpha^2+\gamma^2&0&0&0\\
		0&0&\alpha^2-\beta^2&-\alpha\gamma&0\\
		0&0&-\alpha\gamma&-\gamma^2&0\\
		0&0&0&0&-\beta^2
		\end{pmatrix}}. \label{eqn8}
	\end{equation}
\end{itemize}
\subsubsection{Case of $\mathfrak{g}=A_{5,6}$}
In the basis $\mathcal{B}_5=\{v_1,v_2,v_3,v_4,v_5\}$ in proposition \ref{prpo7} a direct computation using \eqref{riccicur} gives
\begin{equation}
[Ric_{\eta\langle\cdotp,\cdotp\rangle}]=[Ric_{\langle\cdotp,\cdotp\rangle}]=-\frac{1}{2}{\renewcommand{\arraystretch}{0.4}\begin{pmatrix}
	\alpha^2+\beta^2+\gamma^2+\delta^2+\varepsilon^2&\delta\sigma&0&0&0\\
	\delta\sigma&\alpha^2+\beta^2&\beta\gamma&0&0\\
	0&\beta\gamma&\gamma^2+\delta^2+\sigma^2-\alpha^2&
	\delta\varepsilon-\alpha\beta&0\\ 
	0&0&\delta\varepsilon-\alpha\beta&\varepsilon^2-\beta^2-\gamma^2&-\delta\gamma\\
	0&0&0&-\delta\gamma&-\delta^2-\varepsilon^2-\sigma^2
	\end{pmatrix}}. \label{eqn13}
\end{equation}

\subsubsection{Case of $\mathfrak{g}=A_{5,5}$}
In the basis $\mathcal{B}_6=\{v_1,v_2,v_3,v_4,v_5\}$ in proposition \ref{prpo8} a direct computation using \eqref{riccicur} gives
\begin{equation}
[Ric_{\eta\langle\cdotp,\cdotp\rangle}]=[Ric_{\langle\cdotp,\cdotp\rangle}]=-\frac{1}{2}{\renewcommand{\arraystretch}{0.4}\begin{pmatrix}
	\alpha^2+\beta^2+\gamma^2&\gamma\delta&-\beta\delta&-\beta\varepsilon&0\\
	\gamma\delta&\alpha^2+\beta^2+\delta^2+\varepsilon^2&\beta\gamma&0&0\\
	-\beta\delta&\beta\gamma&\gamma^2+\delta^2&\delta\varepsilon&0\\
	-\beta\varepsilon&0&\delta\varepsilon&\varepsilon^2-\alpha^2&-\alpha\beta\\
	0&0&0&-\alpha\beta&-\beta^2-\gamma^2-\delta^2-\varepsilon^2
	\end{pmatrix}}. \label{eqn18}
\end{equation} 

\subsubsection{Case of $\mathfrak{g}=A_{5,3}$}
In the basis $\mathcal{B}_7=\{v_1,v_2,v_3,v_4,v_5\}$ in proposition \ref{prpo08}, a direct computation using \eqref{riccicur} gives
\begin{equation}
[Ric_{\eta\langle\cdotp,\cdotp\rangle}]=[Ric_{\langle\cdotp,\cdotp\rangle}]=-\frac{1}{2}{\renewcommand{\arraystretch}{0.4}\begin{pmatrix}
	\alpha^2+\beta^2+\gamma^2+\delta^2&\delta\varepsilon&0&0&0\\
	\delta\varepsilon&\alpha^2+\beta^2+\varepsilon^2&\beta\gamma&0&0\\
	0&\beta\gamma&\gamma^2+\delta^2+\varepsilon^2-\alpha^2&-\alpha\beta&0\\
	0&0&-\alpha\beta&-\beta^2-\gamma^2&-\delta\gamma\\
	0&0&0&-\delta\gamma&-\delta^2-\varepsilon^2
	\end{pmatrix}}. \label{eqn18}
\end{equation}

\subsubsection{Case of $\mathfrak{g}=A_{5,1}$}
In the basis $\mathcal{B}_8=\{v_1,v_2,v_3,v_4,v_5\}$ in proposition \ref*{prpo008}, a direct computation using \eqref{riccicur} gives
\begin{equation}
[Ric_{\eta\langle\cdotp,\cdotp\rangle}]=[Ric_{\langle\cdotp,\cdotp\rangle}]=-\frac{1}{2}{\renewcommand{\arraystretch}{0.4}\begin{pmatrix}
	\alpha^2+\beta^2+\gamma^2&0&0&0&0\\
	0&\alpha^2+\beta^2&\beta\gamma&0&0\\
	0&\beta\gamma&\gamma^2&0&0\\
	0&0&0&-\alpha^2&-\alpha\beta\\
	0&0&0&-\alpha\beta&-\beta^2-\gamma^2
	\end{pmatrix}}. \label{eqn0018}
\end{equation}

\subsubsection{Case of $\mathfrak{g}=A_{5,2}$}
In the basis $\mathcal{B}_9=\{v_1,v_2,v_3,v_4,v_5\}$ in proposition \ref*{prpo8}, a direct computation using \eqref{riccicur} gives
\begin{equation}
[Ric_{\eta\langle\cdotp,\cdotp\rangle}]=[Ric_{\langle\cdotp,\cdotp\rangle}]=-\frac{1}{2}{\renewcommand{\arraystretch}{0.4}\begin{pmatrix}
	\alpha^2+\beta^2+\gamma^2+\delta^2&0&0&0&0\\
	0&\alpha^2+\beta^2&\beta\gamma&0&0\\
	0&\beta\gamma&\gamma^2-\alpha^2&-\alpha\beta&0\\
	0&0&-\alpha\beta&\delta^2-\beta^2-\gamma^2&0\\
	0&0&0&0&-\delta^2
	\end{pmatrix}}. \label{eqn23}
\end{equation} 

\section{Ricci curvature equation on $5$-dimensional nilpotent Lie groups}
In this section, we study the existence of solution of \eqref{riccic} on any $5$-dimensional nilpotent Lie group. The theorems obtained arise from the study of equivalents polynomial systems that we established from the results of the previous section. We give the complete proofs of theorems \ref{th1} and \ref{th7}. By applying the same technique on the other algebras with the equivalent polynomial systems we obtain the proofs of theorems \ref{th5}, \ref{th6}, \ref{th8}, \ref{th9}, \ref{th10} and \ref{th11} .   \\
In the following $G$ denote the nilpotent Lie group of Lie algebra $\mathfrak{g}$ where $\mathfrak{g}$ represent each of the $5$-dimensional Lie algebras.  
\subsection{Case of the Lie algebra  $\mathfrak{g}=5A_{1}$}
\begin{theorem}
	Let T be a left-invariant symmetric tensor field of type (0,2) on $G$. there exists a pair $(g,t)$ where $g$ is a left-invariant Riemanian metric on $G$ and t is a non-zero real constant such that $Ric(g)= t^2T$ if and only if $T=0$. \label{th3}
	
\end{theorem}
\begin{prof}
	The form of $T$ in the Basis $\mathcal{B}_1$ of $5A_{1}$ is impose by that of the matrix expression of the Ricci curvature. As long as $T$ and the Ricci curvature are of the same form in B,\eqref{riccic} is equivalent to $T=0.$ 
\end{prof}
\subsection{Case of the Lie algebra  $\mathfrak{g}=A_{5,4}$}
\begin{theorem}
	Let T be a left-invariant symmetric tensor field of type (0,2) on $G$. there exists a pair $(g,t)$ where $g$ is a left-invariant Riemanian metric on $G$ and t is a non-zero real constant such that  $Ric(g)= t^2T$ if and only if T is of the following form in a basis $\mathcal{B}_2=\{v_1,v_2,v_3,v_4,v_5\}$ described in proposition \ref{prop00001}:
	\begin{equation*}
	T={\renewcommand{\arraystretch}{0.4}\begin{pmatrix}
		a&f&0&0&0\\
		f&b&0&0&0\\
		0&0&c&l&0\\
		0&0&l&d&0\\
		0&0&0&0&e
		\end{pmatrix}}.
	\end{equation*} with the following conditions:
	\begin{enumerate}
		\item[(1)]  $a+b+e=0$
		\quad(2) $b<0$
		\quad(3) $d<0$
		\quad(4) $b-c\geq0$.
		\item[(5)] $f\pm \sqrt{-b(b-c)}=0$\quad
		(6) $l\pm \sqrt{-d(b-c)}=0$.
		
	\end{enumerate}\label{th1}
\end{theorem} 
\begin{prof}
	The matrix expression of the Ricci curvature in the Basis $\mathcal{B}_2$ of $A_{5,4}$ induce the form of $T$. Since $T$ and the Ricci curvature are of the same form in $\mathcal{B}_2 $, \eqref{riccic} is equivalent to the following algebraic system.
	\begin{eqnarray} 
	\alpha^2+\beta^2+2at^2&=&0\nonumber\\
	\gamma^2+2bt^2&=&0\nonumber\\
	\alpha^2+\gamma^2+2ct^2&=&0\nonumber\\
	\beta^2+2dt^2&=&0\nonumber\\	\alpha^2+\beta^2+\gamma^2-2et^2&=&0\label{systemaA5,4}\\
	\alpha\gamma+2ft^2&=&0\nonumber\\
	\alpha\beta+2lt^2&=&0\nonumber\\
	\alpha\in \mathbb{R}\; \text{and} \;\beta,\gamma&>&0\nonumber
	\end{eqnarray}
The polynomial system made up of the first five equations of system  \eqref{systemaA5,4} is linear in $\alpha^2,\beta^2,\gamma^2$ and $t^2$.
It's a system of five equations relating in four unknowns. Therefore there will appear $5-4=1$ condition of compatibility of the system to have solutions.The extended matrix of the system is the following:

$${\renewcommand{\arraystretch}{0.4}\begin{pmatrix}
	1&1&1&-2e&\bigm|&0\\
	0&1&0&2d&\bigm|&0\\
	0&0&1&2b&\bigm|&0\\
	1&1&0&2a&\bigm|&0\\
	1&0&1&2c&\bigm|&0
	\end{pmatrix}} $$	
The Gauss algorithm applied to the extended matrix gives the following matrix:
$${\renewcommand{\arraystretch}{0.4}\begin{pmatrix}
	1&1&1&-2e&\bigm|&0\\
	0&1&0&2d&\bigm|&0\\
	0&0&1&2b&\bigm|&0\\
	0&0&0&2(a+b+e)&\bigm|&0\\
	0&0&0&0&\bigm|&0
	\end{pmatrix}} $$
therefore The polynomial system made up of the first five equations of system  \eqref{systemaA5,4} admits a solution if and only if the relationship $a+b+e=0$ is satisfied. \\
Under the condition $a+b+e=0$, the set of solutions of the previous polynomial system is the set of vectors of the form
$${\renewcommand{\arraystretch}{0.4}\begin{pmatrix}
	\alpha^2\\
	\beta^2\\
	\gamma^2\\
	t^2
	\end{pmatrix}} ={\renewcommand{\arraystretch}{0.4}\begin{pmatrix}
	2(b+d+e)t^2\\
	-2bt^2\\
	-2dt^2\\
	t^2
	\end{pmatrix}}\quad with \; t\in\mathbb{R}^*$$.	

Using the condions $\alpha\in \mathbb{R}$, $\beta,\gamma>0$ and the sixth and seventh equations of the system \eqref{systemaA5,4} we obtain that for every $t\in \mathbb{R}^*$, there exist $\alpha,\beta$ and $\gamma$ dependent of $t$ such that $(\alpha,\beta,\gamma,t)$ is a solution of \eqref{systemaA5,4} if and only if the conditions  $(1),(2),(3),(4),(5)$ and $(6)$ of the theorem hold.  

\end{prof}

\subsection{Case of the Lie algebra  $\mathfrak{g}=A_{3,1}\oplus 2A_1$}
\begin{theorem}
	Let T be a left-invariant symmetric tensor field of type (0,2) on $G$. there exists a pair $(g,t)$ where $g$ is a left-invariant Riemanian metric on $G$ and t is a non-zero real constant such that Ric  $(g)= t^2T$ if and only if T is of the following form in a basis $\mathcal{B}_3=\{v_1,v_2,v_3,v_4,v_5\}$ described in the proposition \ref{prop00002}
	\begin{equation*}
	T=diag(a,b,0,0,c).
	\end{equation*} with the following conditions:
	\begin{itemize}
		\item [(1)] $a<0$\quad	(2) $a=b=-c$.
	\end{itemize}\label{th5}
\end{theorem}
\begin{prof}
	The matrix expression of the Ricci curvature in the Basis $\mathcal{B}_3$ of $A_{3,1}\oplus 2A_1$ induce the form of $T$. Since $T$ and the Ricci curvature are of the same form in $\mathcal{B}_3$, \eqref{riccic} is equivalent to the following algebraic system.
	\begin{eqnarray} 
	\alpha^2+2at^2&=&0\nonumber\\
	\alpha^2+2bt^2&=&0\label{system2}\\
	\alpha^2-2ct^2&=&0\nonumber\\
	\alpha&>&0\nonumber
	\end{eqnarray}
	From the system below, $\alpha^2=2ct^2=-2bt^2=-2at^2$. Then whith the condition $\alpha>0$, we obtain that for every $t\in \mathbb{R}^*$, there exist $\alpha$ dependent of $t$ such that $(\alpha,t)$ is solution of \eqref{system2} if and only if the conditions $(1)$ and $(2)$ of the theorem hold.
\end{prof}
\subsection{Case of the Lie algebra  $\mathfrak{g}=A_{4,1}\oplus A_1$}
\begin{theorem}
	Let T be a left-invariant symmetric tensor field of type (0,2) on $G$. there exists a pair $(g,t)$ where $g$ is a left-invariant Riemanian metric on $G$ and t is a non-zero real constant such that $Ric(g)= t^2T$ if and only if $T$ is of the form  $T={\renewcommand{\arraystretch}{0.4}\begin{pmatrix}
		a&0&0&0&0\\
		0&b&e&0&0\\
		0&e&c&0&f\\
		0&0&0&0&0\\
		0&0&f&0&d
		\end{pmatrix}} $ in a  basis $\mathcal{B}_4=\{v_1,v_2,v_3,v_4,v_5\}$ described in proposition \ref{prpo6} with the following conditions:
	\begin{itemize}
		\item [(1)]$b+c+d=0$\quad(2) $A=a-2b-c\geq0$\quad(3) $B=b-a>0$\quad(4) $C=a+d<0$.
		\item [(5)] $e\pm \sqrt{AB}=0$\quad(6)$f\pm \sqrt{-AC}=0$.	
	\end{itemize}
	or  $T$ is of the form  $T={\renewcommand{\arraystretch}{0.4}\begin{pmatrix}
		a&0&0&0&0\\
		0&b&e&0&0\\
		0&e&c&0&f\\
		0&0&0&0&0\\
		0&0&f&0&d
		\end{pmatrix}} $ in a basis $\mathcal{B}_4=\{v_1,v_2,v_3,v_4,v_5\}$ described in proposition \ref{prpo6} with the following conditions:
	\begin{itemize}
		\item[(1)] $2b-c-a+d=0$\quad(2) $b-a-e=0$ \quad(3) $d\geq0$\quad(4) $e>0$\quad (5) $b+d<0$\quad(6)$f\pm \sqrt{-d(b+d)}=0$.	
	\end{itemize}
	\label{th6}

\end{theorem}
\begin{prof}	
	Since $T$ and the Ricci curvature are of the same form in $\mathcal{B}_4$, \eqref{riccic} is equivalent to the algebraic system \eqref{systemA5,40} in the first case and the system \eqref{systemA5,41} in the second case. The proof of this theorem is similar to the one of theorem \ref{th1}	
	 \begin{itemize}
		\item[1.]\underline{First case}
	\begin{align}
	\alpha^2+\beta^2+\gamma^2+2at^2&=0\nonumber\\
	\alpha^2+\gamma^2+2bt^2&=0\nonumber\\
	\beta^2-\alpha^2+2ct^2&=0\nonumber\\
	\beta^2+\gamma^2-2dt^2&=0\label{systemA5,40}\\
	\beta\gamma+2et^2&=0\nonumber\\
	\alpha\gamma-2ft^2&=0\nonumber\\
	\alpha,\beta>0\; \text{and} \;\gamma\in \mathbb{R}\nonumber
	\end{align}
	\item[2.]\underline{Second case}
	\begin{align}
	\alpha^2+\beta^2+\gamma^2+2at^2&=0\nonumber\\
	\alpha^2+\gamma^2+2bt^2&=0\nonumber\\
	-\beta^2+\alpha^2+2ct^2&=0\nonumber\\
	\gamma^2-2dt^2&=0\label{systemA5,41}\\
	\beta^2-2et^2&=0\nonumber\\
	\alpha\gamma-2ft^2&=0\nonumber\\
	\alpha,\beta>0\; \text{and} \;\gamma\in \mathbb{R}\nonumber
	\end{align}
\end{itemize}
%
%
%
\end{prof} 

\subsection{Case of the Lie algebra  $\mathfrak{g}=A_{5,6}$}
\begin{theorem}
	Let T be a left-invariant symmetric tensor field of type (0,2) on $G$. If T is of the form $T={\renewcommand{\arraystretch}{0.4}\begin{pmatrix}
		a&f&0&0&0\\
		f&b&g&0&0\\
		0&g&c&h&0\\
		0&0&h&d&i\\
		0&0&0&i&e
		\end{pmatrix}}$ in a basis $\mathcal{B}_5=\{v_1,v_2,v_3,v_4,v_5\}$ described in proposition \ref{prpo7}
	and
	\begin{itemize}
		\item [(1)] $b+c+d+e=0$\quad(2) $i^2-f^2\neq0$\quad(3) $A+D-E>0$\quad(4) $B-D+E>0$\quad(5) $C+D>0$.
		\item[(6)] $e-D-E\geq0$\quad(7) $D>0$\quad(8) $E>0$\quad(9) $BC-DC+CE+BD+ED-D^2-g^2=0$.
		\item[(10)] $l\pm\sqrt{(B-D+E)(C+D)}=0$\quad(11) $f\pm\sqrt{(e-D-E)D}=0$.
		\item[(12)] $i\pm\sqrt{(e-D-E)(C+D)}=0$\quad(13) $h\pm\sqrt{(e-D-E)E}\pm\sqrt{(A+D-E)(B-D+E)}=0$.	
		\\ with $A=b+c-a$, $B=a-2b-c$, $C=2b+c+d-a$.
		$D=\frac{f^2C}{i^2-f^2}$ and $E=e-\frac{f^2C}{i^2-f^2}-\frac{i^2-f^2}{C}$.		
	\end{itemize}	
	then there exists a pair $(g,t)$ where $g$ is a left-invariant Riemanian metric on $G$ and t is a non-zero real constant such that Ric  $(g)= t^2T$.\label{th7}
\end{theorem}
\begin{prof}
	The form of $T$ in the Basis $\mathcal{B}_5$ of $A_{5,6}$ is impose by that of the matrix expression of the Ricci curvature. As long as $T$ and the Ricci curvature are of the same form in B, \eqref{riccic} is equivalent to the following algebraic system.
	\begin{eqnarray} 
		\alpha^2+\beta^2+\gamma^2+\delta^2+\varepsilon^2+2at^2&=&0\nonumber\\
	\alpha^2+\beta^2+2bt^2&=&0\nonumber\\
	\gamma^2+\delta^2+\sigma^2-\alpha^2+2ct^2&=&0\nonumber\\
	\varepsilon^2-\beta^2-\gamma^2+2dt^2&=&0\nonumber\\
	\delta^2+\varepsilon^2+\sigma^2-2et^2&=&0\label{system5.6}\\
	\delta\sigma+2ft^2&=&0\nonumber\\
	\beta\gamma+2lt^2&=&0\nonumber\\
	\delta\varepsilon-\alpha\beta+2ht^2&=&0\nonumber\\ 
	\delta\gamma-2it^2&=&0\nonumber\\
	\alpha<0;\;\gamma,\varepsilon,\sigma>0\; \text{and} \;\beta,\delta&\in&\mathbb{R}\nonumber
	\end{eqnarray}
	
%
	Under the conditions $(1),(2),(3),(4),(5),(6),(7),(8),(9),(10),(11),(12)$ and $(13)$ a set of vectors of the following form  is a set of solutions of system \eqref{system5.6}.
	$${\renewcommand{\arraystretch}{1}\begin{pmatrix}
		\alpha\\
		\beta\\
		\gamma\\
		\delta\\
		\varepsilon\\
		\sigma\\
		t
		\end{pmatrix}} ={\renewcommand{\arraystretch}{1}\begin{pmatrix}
		-t\sqrt{2(A+D-E)}\\
		\pm t\sqrt{2(B+E-D)}\\
		t\sqrt{2(C+D)}\\
		\pm t\sqrt{2(e-D-E)}\\
	t\sqrt{2E}\\
		t\sqrt{2D}\\
		t
		\end{pmatrix}} \;\forall\;t\in \mathbb{R}^*.$$	
	
%
%
\end{prof}
\subsection{Case of the Lie algebra  $\mathfrak{g}=A_{5,5}$}
\begin{theorem}
	Let T be a left-invariant symmetric tensor field of type (0,2) on $G$. there exists a pair $(g,t)$ where $g$ is a left-invariant Riemanian metric on $G$ and t is a non-zero real constant such that Ric  $(g)= t^2T$ if and only if T is of the following form in a basis $\mathcal{B}_6=\{v_1,v_2,v_3,v_4,v_5\}$ described in the proposition \ref{prpo8}:
	\begin{equation*}
	T={\renewcommand{\arraystretch}{0.4}\begin{pmatrix}
		a&f&l&h&0\\
		f&b&i&0&0\\
		l&i&c&j&0\\
		h&0&j&d&k\\
		0&0&0&k&e
		\end{pmatrix}}.
	\end{equation*} with the following conditions:
	\begin{enumerate}[(i)]
		\item[(1)]  $A=-(a+b+c+2d+2e)>0$
		\quad(2) $B=a+d+e\geq0$
		\quad(3) $C=-(a+c+d+e)>0$.
		\item[(4)] $D=a+2c+3d+3e\geq0$
		\quad(5) $E=-(a+c+2d+2e)>0$
		\quad(6) $f\pm \sqrt{BC}=0$
		\quad(7) $l\pm \sqrt{BD}=0$.
		\item[(8)] $h\pm \sqrt{AD}=0$
		\quad(9) $i\pm \sqrt{DC}=0$
		\quad(10) $j\pm \sqrt{AB}=0$
		\quad(11) $k\pm \sqrt{DE}=0$.	
	\end{enumerate}
\label{th8}
\end{theorem}
\begin{prof}
	%
	On $A_{5,5}$ \eqref{riccic} is equivalent to the following system in a basis $\mathcal{B}_6$ described in the proposition \ref{prop00001}
	\begin{eqnarray} 
	\alpha^2+\beta^2+\gamma^2+2at^2&=&0\nonumber\\
	\alpha^2+\beta^2+\delta^2+\varepsilon^2+2bt^2&=&0\nonumber\\
	\gamma^2+\delta^2+2ct^2&=&0\nonumber\\
	\varepsilon^2-\alpha^2+2dt^2&=&0\nonumber\\
	\beta^2+\gamma^2+\delta^2+\varepsilon^2-2et^2&=&0\nonumber\\
	\gamma\delta+2ft^2&=&0\label{system6}\\
	\beta\delta-2lt^2&=&0\nonumber\\
	\beta\varepsilon-2ht^2&=&0\nonumber\\
	\beta\gamma+2it^2&=&0\nonumber\\
	\delta\varepsilon+2jt^2&=&0\nonumber\\
	\alpha\beta-2kt^2&=&0\nonumber\\
	\alpha,\gamma,\varepsilon>0\; \text{and} \;\beta,\delta&\in& \mathbb{R}\nonumber
	\end{eqnarray}
The proof is similar to that of theorem \ref{th1} 
\end{prof}
\subsection{Case of the Lie algebra  $\mathfrak{g}=A_{5,3}$}
\begin{theorem}
	Let T be a left-invariant symmetric tensor field of type (0,2) on $G$. there exists a pair $(g,t)$ where $g$ is a left-invariant Riemanian metric on $G$ and t is a non-zero real constant such that Ric  $(g)= t^2T$ if and only if T is of the following form in a basis $\mathcal{B}_7=\{v_1,v_2,v_3,v_4,v_5\}$ described in the proposition \ref{prpo08}
	\begin{equation*}
	T={\renewcommand{\arraystretch}{0.4}\begin{pmatrix}
		a&f&0&0&0\\
		f&b&l&0&0\\
		0&l&c&h&0\\
		0&0&h&d&i\\
		0&0&0&i&e
		\end{pmatrix}}.
	\end{equation*} with the following conditions:
	\begin{enumerate}[(i)]
		\item[(1)] $A=-(b+c+d+e)>0$
		\quad(2) $B=b+c+d+2e\geq0$
		\quad(3) $C=-(a+b+2c+2d+3e)>0$.
		\item[(4)] $D=a+b+2c+3d+3e\geq0$
		\quad(5) $E=-(a+b+c+2d+2e)>0$
		\quad(6) $f\pm \sqrt{AB}=0$.
		\item[(7)] $l\pm \sqrt{BD}=0$
		\quad(8) $-h\pm \sqrt{DE}=0$
		\quad(9) $-i\pm \sqrt{BC}=0$.	
	\end{enumerate}\label{th9}
\end{theorem}
\begin{prof}
	In a basis $\mathcal{B}_7$ described in the proposition \ref{prpo08}, \eqref{riccic} is equivalent to the following algebraic system.
	\begin{eqnarray} 
	\alpha^2+\beta^2+\gamma^2+\delta^2+2at^2&=&0\nonumber\\
	\alpha^2+\beta^2+\varepsilon^2+2bt^2&=&0\nonumber\\
	\gamma^2+\delta^2+\varepsilon^2-\alpha^2+2ct^2&=&0\nonumber\\
	\beta^2+\gamma^2-2dt^2&=&0\nonumber\\
	\delta^2+\varepsilon^2-2et^2&=&0\nonumber\\
	\delta\varepsilon+2ft^2&=&0\label{system8}\\
	\beta\gamma+2gt^2&=&0\nonumber\\
	\alpha\beta-2ht^2&=&0\nonumber\\
	\delta\gamma-2it^2&=&0\nonumber\\
	\alpha,\gamma,\varepsilon>0\; \text{and} \;\beta,\delta&\in& \mathbb{R}\nonumber
	\end{eqnarray}
	The proof is similar to that of theorem \ref{th1} 
\end{prof}

\subsection{Case of the Lie algebra  $\mathfrak{g}=A_{5,1}$}
\begin{theorem}
	Let T be a left-invariant symmetric tensor field of type (0,2) on $G$. there exists a pair $(g,t)$ where $g$ is a left-invariant Riemanian metric on $G$ and t is a non-zero real constant such that Ric  $(g)= t^2T$ if and only if T is of the following form in a basis $\mathcal{B}_8=\{v_1,v_2,v_3,v_4,v_5\}$ described in the proposition \ref{prpo008}
	\begin{equation*}
	T={\renewcommand{\arraystretch}{0.4}\begin{pmatrix}
		a&0&0&0&0\\
		0&b&f&0&0\\
		0&f&c&0&0\\
		0&0&0&d&l\\
		0&0&0&l&e
		\end{pmatrix}}.
	\end{equation*}and the following conditions are satisfied:
	\begin{enumerate}[(i)]
		\item[(1)]  $a+b+c=0$
		\quad(2) $d>0$
		\quad(3) $c<0$
		\quad(4) $b+d\leq0$.
		\item[(5)] $f\pm \sqrt{c(b+d)}=0$
		\quad(6) $l\pm \sqrt{-d(b+d)}=0$.	
	\end{enumerate}\label{th10}
\end{theorem}
\begin{prof}
	In a basis $\mathcal{B}_8$ described in the proposition \ref{prpo008}, \eqref{riccic} is equivalent to the following algebraic system.
	\begin{eqnarray} 
		\alpha^2+\beta^2+\gamma^2+2at^2&=&0\nonumber\\
	\alpha^2+\beta^2+2bt^2&=&0\nonumber\\
		\gamma^2+2ct^2&=&0\nonumber\\
		\alpha^2-2dt^2&=&0\label{system8}\\
	\beta^2+\gamma^2-2et^2&=&0\nonumber\\	
	\beta\gamma+2ft^2&=&0\nonumber\\
	\alpha\beta-2lt^2&=&0\nonumber\\
	\alpha,\gamma>0\; \text{and} \;\beta&\in& \mathbb{R}\nonumber
	\end{eqnarray}
	The proof is similar to that of theorem \ref{th1}
\end{prof}

\subsection{Case of the Lie algebra  $\mathfrak{g}=A_{5,2}$}
\begin{theorem}
	Let T be a left-invariant symmetric tensor field of type (0,2) on $G$. there exists a pair $(g,t)$ where $g$ is a left-invariant Riemanian metric on $G$ and t is a non-zero real constant such that Ric  $(g)= t^2T$ if and only if T is of the following form in a basis $\mathcal{B}_9=\{v_1,v_2,v_3,v_4,v_5\}$ described in the proposition \ref{prpo10}
	\begin{equation*}
	T={\renewcommand{\arraystretch}{0.4}\begin{pmatrix}
		a&0&0&0&0\\
		0&b&f&0&0\\
		0&f&c&l&0\\
		0&0&l&d&0\\
		0&0&0&0&e
		\end{pmatrix}}.
	\end{equation*}and the following conditions are satisfied:
	\begin{enumerate}
		\item[(1)]  $b+c+d+e=0$\quad (2) $e>0$
		\quad (3) $A=a-b+e<0$\quad
		(4) $B=a-b+d+2e\geq0$.
		\item[(5)] $C=a+d+2e<0$\quad
		(6) $f\pm \sqrt{-AB}=0$\quad
		(7) $l\pm \sqrt{-BC}=0$.	
	\end{enumerate}
\label{th11}
\end{theorem}
\begin{prof}
	in a basis $\mathcal{B}_9$ described in the proposition \ref{prpo10}, \eqref{riccic} is equivalent to the following algebraic system.
	\begin{eqnarray} 
	\alpha^2+\beta^2+\gamma^2+\delta^2+2at^2&=&0\nonumber\\
	\alpha^2+\beta^2+2bt^2&=&0\nonumber\\
	\gamma^2-\alpha^2+2ct^2&=&0\label{system9}\\
	\delta^2-\beta^2-\gamma^2+2dt^2&=&0\nonumber\\
		\delta^2-2et^2&=&0\nonumber\\
			\beta\gamma+2ft^2&=&0\nonumber\\
	\alpha\beta-2lt^2&=&0\nonumber\\
	\alpha,\gamma,\delta>0\;\text{and}\;\beta&\in& \mathbb{R}\nonumber
	\end{eqnarray}
	The proof is similar to that of theorem \ref{th1}
\end{prof}

\end{document}